\documentclass[preprint,12pt]{elsarticle}
\usepackage{geometry}
\newif\ifsubmissionpackage
\submissionpackagetrue % set to \submissionpackagefalse to obtain the clean manuscript for upload
\usepackage{lmodern}
\usepackage{anyfontsize}
\usepackage[utf8]{inputenc}
\usepackage{amsmath,amssymb,amsthm}
\usepackage{mathrsfs}
\usepackage[shortlabels]{enumitem}
\usepackage{booktabs}
\usepackage{array}
\usepackage{microtype}
\usepackage[colorlinks=true,linkcolor=blue,citecolor=blue,urlcolor=blue]{hyperref}
\usepackage[capitalise,nameinlink]{cleveref}
\usepackage{float}
\numberwithin{equation}{section}

\theoremstyle{plain}
\newtheorem{theorem}{Theorem}[section]
\newtheorem{proposition}[theorem]{Proposition}
\newtheorem{lemma}[theorem]{Lemma}
\newtheorem{corollary}[theorem]{Corollary}
\theoremstyle{definition}
\newtheorem{definition}[theorem]{Definition}
\newtheorem{example}[theorem]{Example}
\theoremstyle{remark}
\newtheorem{remark}[theorem]{Remark}

\crefname{theorem}{Theorem}{Theorems}
\crefname{proposition}{Proposition}{Propositions}
\crefname{lemma}{Lemma}{Lemmas}
\crefname{corollary}{Corollary}{Corollaries}
\crefname{definition}{Definition}{Definitions}
\crefname{example}{Example}{Examples}
\crefname{remark}{Remark}{Remarks}
\crefname{section}{Section}{Sections}
\crefname{subsection}{Section}{Sections}
\crefname{appendix}{Appendix}{Appendices}
\crefname{table}{Table}{Tables}

\allowdisplaybreaks 

\newcommand{\N}{\mathbb{N}}
\newcommand{\Q}{\mathbb{Q}}
\newcommand{\R}{\mathbb{R}}
\newcommand{\cB}{\mathscr{B}}
\newcommand{\cO}{\mathscr{O}}
\newcommand{\cC}{\mathscr{C}}
\newcommand{\cK}{\mathscr{K}}
\newcommand{\cM}{\mathscr{M}}
\newcommand{\cP}{\mathscr{P}}
\newcommand{\bmu}{\overline{\mu}}
\newcommand{\bnu}{\overline{\nu}}
\newcommand{\bv}{\overline{v}}
\newcommand{\bla}{\overline{\lambda}}
\newcommand{\Lu}{\textup{(L)}}
\newcommand{\LK}{\textup{(L$_K$)}}
\newcommand{\PL}{\textup{(PL)}}
\newcommand{\PLK}{\textup{(PL$_K$)}}
\newcommand{\CBX}{\textup{(CBX)}}
\newcommand{\CS}{\textup{(CS)}}
\newcommand{\DW}{\textup{(DW)}}
\newcommand{\DP}{\textup{(DP)}}
\newcommand{\DS}{\textup{(DS)}}
\newcommand{\DE}{\textup{(DE)}}
\newcommand{\PR}{\textup{(PR)}}
\newcommand{\1}{\mathbf{1}}
\newcommand{\core}{\operatorname{core}}
\newcommand{\cl}{\operatorname{cl}}
\newcommand{\inter}{\operatorname{int}}
\newcommand{\dist}{\operatorname{dist}}
\newcommand{\upto}{\uparrow}
\newcommand{\downto}{\downarrow}

\begin{document}

\begin{frontmatter}

\title{Lusin and pseudo-Lusin theorems for monotone measures: necessary conditions, conjugate duality and a correction of a theorem of Castaldo and Marinacci}

\author[a,b]{Duc Toan Truong}
\ead{toan.truongduc@vlu.edu.vn}

\author[b,c]{Kiet Ho}
\ead{hoanhkiet@vlu.edu.vn}

\author[d]{Thanh Ngo Vo Ke}
\ead{nvk.thanh@hutech.edu.vn}

\author[e]{Cuong Le Quoc}
\ead{cuonglq@ptit.edu.vn}

\author[f]{Tran Nhat Luan}
\ead{trannhatluanphpy@gmail.com, luan.tn@vinuni.edu.vn}

\address[a]{Laboratory for Chemical Computation and Modeling, Institute for Computational Science and Artificial Intelligence, Van Lang University, Ho Chi Minh City 70000, Vietnam}

\address[b]{Atomic Molecular and Optical Physics Research Group, Science and Technology Advanced Institute, Van Lang University, Ho Chi Minh City, 70000, Vietnam}

\address[c]{Faculty of Applied Technology, Van Lang School of Technology, Van Lang University, Ho Chi Minh City 70000, Vietnam}

\address[d]{Institute of engineering, HUTECH University, Ho Chi Minh City, Vietnam}

\address[e]{Posts and Telecommunications Institute of Technology - Ho Chi Minh City Campus, Ho Chi Minh City 70000, Vietnam}

\address[f]{Center for Environmental Intelligence, VinUniversity, Hanoi 10000, Vietnam}

\fntext[1]{Corresponding authors: Thanh Ngo Vo Ke, Tran Nhat Luan}

\begin{abstract}
Lusin's theorem for monotone measures has been studied mainly through sufficient conditions. Here we study it from the side of necessary conditions and of conjugate duality. We show that, for an arbitrary monotone measure on a metric space, the Lusin property with closed sets implies regularity, and the Lusin property with compact sets implies order continuity and tightness. Consequently, for finite submeasures the Lusin property is characterized by order continuity on compact spaces, and the compact-set Lusin property is characterized by order continuity on Polish spaces. We then consider a pseudo-Lusin property, in which the Lusin set is required to have almost full measure, and observe that it is exactly the Lusin property of the conjugate measure. For finite monotone measures whose conjugate is subadditive, a class including supermodular capacities, belief functions and coherent lower probabilities, the pseudo-Lusin property is equivalent to continuity from below on compact spaces, and so is its compact-set version on Polish spaces. As an application, we give a corrected form, together with its converse, of a Lusin-type theorem of Castaldo and Marinacci whose original statement was refuted by Wiesel. Several examples, one of them based on meagre sets, show that the hypotheses cannot be omitted.
\end{abstract}


\end{frontmatter}

\tableofcontents

\section{Introduction}\label{sec:intro}

Lusin's theorem states that a Borel function on a metric space equipped with a finite Borel measure is continuous on a closed set whose complement has arbitrarily small measure. Since the pioneering work of Wu and Ha \cite{WH94} on the regularity of fuzzy measures on metric spaces, the theorem and the closely related regularity properties have been extended to many classes of non-additive set functions. Jiang and Suzuki \cite{JS96} proved regularity and Lusin's theorem for exhaustive fuzzy measures that are autocontinuous from above, and, on complete separable metric spaces, tightness and a version with compact Lusin sets; Jiang, Wang and Ziou \cite{JWZ99} obtained the corresponding results for lower semicontinuous fuzzy measures with the pseudometric generating property (p.g.p.). Li and Yasuda \cite{LY04} proved Lusin's theorem for finite continuous weakly null-additive fuzzy measures, and Li and Mesiar \cite{LM11} established it for monotone measures fulfilling condition (E) and having the p.g.p., among other sufficient conditions which are listed in \cite[Section~7]{LM11}. Further results concern the regularity of non-additive measures on complete separable metric spaces \cite{Cao15,WYT13}, measures with values in ordered topological vector spaces \cite{WT14}, and uncertain measures \cite{YZH24}. We refer to the survey \cite{LMPK15} and to the monographs \cite{Den94,Pap95,WK09} for background.

In most of these works the Lusin property is the \emph{goal}, and regularity together with some continuity is the \emph{tool}. It seems natural to ask to what extent these tools are also necessary. For particular classes this question has already been addressed: Wiesel \cite{Wie22} proved that, for normalized sub-additive capacities (continuous from below) on a compact metric space, the compact-set Lusin property is equivalent to continuity from above along decreasing sequences of open sets, and Yang et. al. \cite{YZH24} obtained equivalences between Lusin's theorem and continuity conditions for uncertain measures on compact metric spaces.

A second motivation comes from imprecise probability. Castaldo and Marinacci \cite{CM02} stated a Lusin theorem for supermodular capacities $v$ on Polish spaces, in which the Lusin set $K$ is compact and has \emph{large} capacity, $v(K)\ge 1-\varepsilon$, rather than a complement of small capacity. Wiesel \cite{Wie22} showed by a counterexample that this statement is not correct as it stands, and located the gap in its proof. For non-additive measures the two requirements "$\mu(X\setminus F)$ small" and "$\mu(F)$ close to $\mu(X)$" are genuinely different, and they are related through the conjugate measure $\bmu(A)=\mu(X)-\mu(X\setminus A)$, which also underlies the pseudo-convergence notions studied in \cite{LMPK15,LOM21}.

The aim of this paper is to put these two lines together. Our contributions may be summarized as follows.
\begin{enumerate}[(1)]
\item \emph{Necessary conditions} (\cref{thm:necessity}). For an arbitrary monotone measure $\mu$ on a metric space (possibly infinite and without any continuity assumption), the Lusin property \Lu{} implies regularity, and the compact-set Lusin property \LK{} implies order continuity and tightness. The first implication is elementary and may well be folklore, but we have not found it stated in the literature; the compactness argument for the second one is the one used by Wiesel \cite{Wie22} for sub-additive capacities.
\item \emph{Characterizations for submeasures} (\cref{cor:compact-sub,thm:compact-Lusin}). For a finite submeasure, \Lu{} is equivalent to order continuity (and to continuity) on compact spaces, and \LK{} is equivalent to order continuity on Polish spaces, and to order continuity plus tightness on arbitrary metric spaces. The sufficiency parts are known \cite{JS96,JWZ99,LY04,LM11,Cao15,WYT13}; the necessity parts seem to be new.
\item \emph{Pseudo-Lusin property and conjugate duality} (\cref{sec:duality,sec:pseudo}). We introduce the pseudo-Lusin property \PL{}, show that \PL$[\mu]$ is exactly \Lu$[\bmu]$, give a dictionary between conditions on $\mu$ and on $\bmu$ (\cref{prop:dictionary}), and show that \Lu{} and \PL{} are independent in general.
\item \emph{Conjugate subadditive measures} (\cref{thm:CS}). For a finite monotone measure $v$ satisfying the condition \CS{} of conjugate subadditivity (which holds for supermodular capacities, belief functions and coherent lower probabilities), \PL$[v]$ is equivalent to continuity from below on compact spaces, and \PLK$[v]$ is equivalent to continuity from below on Polish spaces.
\item \emph{The Castaldo--Marinacci theorem} (\cref{sec:CM}). On Polish spaces, the conclusion of \cite[Theorem~1]{CM02} holds for a supermodular (or merely \CS) normalized capacity \emph{if and only if} it is continuous from below (\cref{thm:CM-corrected}). The sufficiency is, in dual form, already contained in \cite[Theorem~6]{JWZ99}; the necessity and the ``if and only if'' formulation are new. The counterexample (\cref{ex:CM-counter}) and the diagnosis of the gap are due to Wiesel \cite{Wie22}; we only make them more explicit.
\item \emph{Examples.} An example based on meagre and comeagre sets (\cref{ex:meagre}) shows that \CS{} cannot be dropped in \cref{thm:CS}, and that a finite, order continuous, supermodular monotone measure which is continuous from above need not have the Lusin property on a compact space.
\end{enumerate}
We also include some elementary transfer results for pairs of monotone measures (\cref{sec:pairs}) and a routine application to approximation by continuous functions (\cref{sec:approx}). \Cref{tab:status} indicates which results are known and which appear to be new.

\begin{table}[H]
\centering
\small
\caption{Overview of the main results and their status.}\label{tab:status}
\begin{tabular}{@{}p{0.23\textwidth}p{0.40\textwidth}p{0.30\textwidth}@{}}
\toprule
Result & Content & Status \\
\midrule
\Cref{prop:dictionary} & duality dictionary; \PL$[\mu]\Leftrightarrow$ \Lu$[\bmu]$ & definitional; the \PL{} viewpoint appears new \\
\Cref{prop:CS} & the condition \CS & \CS$\Leftrightarrow\bv$ subadditive is folklore \\
\Cref{lem:submeasure} & submeasures: order continuity $\Leftrightarrow$ continuity & folklore \\
\Cref{thm:necessity} & \Lu$\Rightarrow$ regular; \LK$\Rightarrow$ order continuous and tight & new ((a) possibly folklore) \\
\Cref{thm:sub-known} & order continuous finite submeasure $\Rightarrow$ \Lu & known \cite{JS96,JWZ99,LY04,LM11} \\
\Cref{cor:compact-sub,thm:compact-Lusin} & characterizations for submeasures & sufficiency known; necessity new \\
\Cref{thm:CS} & pseudo-Lusin theorem under \CS & (a) known by duality; necessity new \\
\Cref{thm:CM-corrected} & corrected Castaldo--Marinacci theorem & sufficiency \cite{JWZ99} in dual form; necessity new \\
\Cref{ex:CM-counter} and \cref{subsec:gap} & counterexample to \cite[Theorem~1]{CM02}; the gap & due to Wiesel \cite{Wie22} \\
\Cref{ex:meagre} & meagre/comeagre example & new in this context, as far as we know \\
\bottomrule
\end{tabular}
\end{table}

The paper is organized as follows. \Cref{sec:prelim} collects definitions. \Cref{sec:duality} is devoted to conjugate duality and to the condition \CS. \Cref{sec:necessity} contains the necessary conditions, and \cref{sec:submeasures} the characterizations for submeasures. Pseudo-Lusin theorems and the main examples are in \cref{sec:pseudo}, and the theorem of Castaldo and Marinacci is discussed in \cref{sec:CM}. \Cref{sec:pairs,sec:approx} contain the transfer results and the application, and \cref{sec:concl} lists some open problems. Self-contained proofs of some known results, and of one pseudo-Lusin theorem, are given in \cref{app:proofs}.

\section{Preliminaries}\label{sec:prelim}

Throughout, $(X,d)$ is a metric space, $\cO$, $\cC$ and $\cK$ denote the families of open, closed and compact subsets of $X$, and $\cB$ is the Borel $\sigma$-algebra. We write $A^c=X\setminus A$. A \emph{Borel function} is a $\cB$-measurable function $f\colon X\to\R$ (finite-valued). A set is \emph{meagre} if it is of first category in $X$, and \emph{comeagre} if its complement is meagre. A \emph{Polish space} is a separable completely metrizable space. We write $A_n\upto A$ (resp.\ $A_n\downto A$) if $(A_n)$ is non-decreasing (resp.\ non-increasing) with union (resp.\ intersection) $A$.

\subsection{Monotone measures}

\begin{definition}\label{def:basic}
A \emph{monotone measure} on $X$ is a set function $\mu\colon\cB\to[0,\infty]$ with $\mu(\emptyset)=0$ and $\mu(A)\le\mu(B)$ whenever $A\subseteq B$. It is \emph{finite} if $\mu(X)<\infty$. For a finite monotone measure $\mu$, the \emph{conjugate} (or \emph{dual}) \emph{measure} is
\[
\bmu(A)=\mu(X)-\mu(A^c),\quad A\in\cB .
\]
Clearly $\bmu$ is a finite monotone measure with $\bmu(X)=\mu(X)$, and $\overline{\bmu}=\mu$.
\end{definition}
\noindent 
Following \cite[Section~2]{LM11} and \cite[Definition~2.3]{LMPK15}, a monotone measure $\mu$ is said to be
\begin{itemize}[leftmargin=1em]
\item \emph{continuous from below} if $A_n\upto A$ implies $\mu(A_n)\to\mu(A)$;
\item \emph{continuous from above} if $A_n\downto A$ and $\mu(A_{n_0})<\infty$ for some $n_0$ imply $\mu(A_n)\to\mu(A)$;
\item \emph{continuous} if it is continuous from below and from above;
\item \emph{order continuous} if $A_n\downto\emptyset$ implies $\mu(A_n)\to0$;
\item \emph{strongly order continuous} if $A_n\downto A$ and $\mu(A)=0$ imply $\mu(A_n)\to0$;
\item \emph{exhaustive} if $\mu(E_n)\to0$ for every sequence $(E_n)$ of pairwise disjoint sets;
\item \emph{autocontinuous from above} (resp.\ \emph{from below}) if $\mu(F_n)\to0$ implies $\mu(E\cup F_n)\to\mu(E)$ (resp.\ $\mu(E\setminus F_n)\to\mu(E)$) for every $E\in\cB$.
\end{itemize}
All sets in these definitions belong to $\cB$. We shall also use the following condition, for which we have not found a standard name in the literature.

\begin{definition}\label{def:CBX}
A finite monotone measure $\mu$ is \emph{continuous from below at $X$}, briefly \CBX, if $A_n\upto X$ implies $\mu(A_n)\to\mu(X)$.
\end{definition}

\begin{remark}\label{rem:CBX-terminology}
By \cref{prop:dictionary}(c) below, \CBX{} for $\mu$ is exactly order continuity of $\bmu$. It should not be confused with the \emph{pseudo-order continuity} of \cite[Definition~5.6]{LMPK15}, which is a different condition, nor with the \emph{strong order pseudo-continuity} of \cite[Definition~2.3]{LOM21}: $A_n\downto A$ and $\mu(A)=0$ imply $\mu(X\setminus A_n)\to\mu(X)$. The latter is a ``mixed'' condition (hypothesis on $\mu$, conclusion on $\bmu$); see \cref{rem:pseudo-comparison}.
\end{remark}
\noindent 
A monotone measure $\mu$ is
\begin{itemize}[leftmargin=1em]
\item \emph{null-additive} if $\mu(E\cup F)=\mu(E)$ whenever $\mu(F)=0$;
\item \emph{weakly null-additive} (w.n.a.) if $\mu(E)=\mu(F)=0$ implies $\mu(E\cup F)=0$;
\item said to have the \emph{pseudometric generating property} (p.g.p.) \cite{DF80} if for every $\varepsilon>0$ there is $\delta>0$ such that $\mu(E)\vee\mu(F)<\delta$ implies $\mu(E\cup F)<\varepsilon$;
\item said to have \emph{property (S)} \cite{Sun94} if every sequence $(A_n)$ with $\mu(A_n)\to0$ has a subsequence $(A_{n_i})$ with $\mu\bigl(\bigcap_{k}\bigcup_{i\ge k}A_{n_i}\bigr)=0$;
\item said to \emph{fulfil condition (E)} \cite{Li03,LM11} if for every double sequence $(E^{(m)}_n)_{n,m\in\N}\subseteq\cB$ such that $E^{(m)}_n\downto E^{(m)}$ as $n\to\infty$ for each $m$ and $\mu\bigl(\bigcup_m E^{(m)}\bigr)=0$, there are increasing sequences $(n_i)$ and $(m_i)$ of natural numbers with $\lim_{k\to\infty}\mu\bigl(\bigcup_{i\ge k}E^{(m_i)}_{n_i}\bigr)=0$. Condition (E) is necessary and sufficient for Egoroff's theorem \cite{Li03}; see also \cite{MUA04}.
\end{itemize}

\begin{definition}\label{def:submeasure}
A \emph{submeasure} is a monotone measure $\mu$ with $\mu(A\cup B)\le\mu(A)+\mu(B)$ for all $A,B\in\cB$.
\end{definition}
\noindent 
This is finite subadditivity only, and no continuity is assumed; this is the usage of Talagrand \cite{Tal08} in connection with Maharam's problem. (In the usage of Dobrakov \cite{DF80}, a submeasure is in addition continuous; we do \emph{not} assume this.)
\noindent 
A finite monotone measure $v$ is \emph{supermodular} (or \emph{$2$-monotone}, \emph{convex}) if $v(A\cup B)+v(A\cap B)\ge v(A)+v(B)$, and \emph{submodular} (or \emph{$2$-alternating}) if the reverse inequality holds, for all $A,B\in\cB$. \emph{Belief functions} (totally monotone capacities) are supermodular. A \emph{coherent lower probability} is the lower envelope of a nonempty family of finitely additive probabilities \cite{Wal91}. The following condition will play a central role.

\begin{definition}\label{def:CS}
A finite monotone measure $v$ is \emph{conjugate subadditive}, briefly \CS, if
\[
v(A)+v(B)\le v(X)+v(A\cap B)\quad\text{for all }A,B\in\cB .
\]
\end{definition}

\subsection{Regularity, tightness and Lusin properties}

\begin{definition}\label{def:regular}
A monotone measure $\mu$ is \emph{regular} (in the sense of Wu and Ha \cite{WH94}) if for every $A\in\cB$ and $\varepsilon>0$ there are $F\in\cC$ and $G\in\cO$ with $F\subseteq A\subseteq G$ and $\mu(G\setminus F)<\varepsilon$. It is \emph{tight} if for every $\varepsilon>0$ there is $K\in\cK$ with $\mu(X\setminus K)<\varepsilon$ (as in \cite{Cao15,JWZ99,WYT13}). A finite $\mu$ is \emph{pseudo-tight} if for every $\varepsilon>0$ there is $K\in\cK$ with $\mu(K)>\mu(X)-\varepsilon$; equivalently, $\bmu$ is tight.
\end{definition}

\begin{definition}\label{def:Lusin}
Let $\mu$ be a monotone measure on $X$.
\begin{enumerate}[(i)]
\item $\mu$ has the \emph{Lusin property} \Lu$[\mu]$ if for every Borel function $f$ and every $\varepsilon>0$ there is $F\in\cC$ with $\mu(X\setminus F)<\varepsilon$ such that $f|_F$ is continuous.
\item \LK$[\mu]$ is the same property with $F\in\cK$.
\item For finite $\mu$, $\mu$ has the \emph{pseudo-Lusin property} \PL$[\mu]$ if for every Borel function $f$ and every $\varepsilon>0$ there is $F\in\cC$ with $\mu(F)>\mu(X)-\varepsilon$ such that $f|_F$ is continuous.
\item \PLK$[\mu]$ is the same property with $F\in\cK$.
\end{enumerate}
\end{definition}

\begin{remark}\label{rem:Lusin-basic}
(i) Compact sets are closed, so \LK{} implies \Lu{} and \PLK{} implies \PL; on a compact space the two versions coincide. (ii) Whether one requires ``$<\varepsilon$'' or ``$\le\varepsilon$'' in these definitions is immaterial.
\end{remark}

\section{Conjugate duality and the condition (CS)}\label{sec:duality}

\begin{proposition}\label{prop:dictionary}
Let $\mu$ be a finite monotone measure and $\nu:=\bmu$.
\begin{enumerate}[(a)]
\item $\nu(X\setminus F)=\mu(X)-\mu(F)$ for every $F\in\cB$. Consequently \PL$[\mu]\Leftrightarrow$ \Lu$[\bmu]$, \PLK$[\mu]\Leftrightarrow$ \LK$[\bmu]$, \Lu$[\mu]\Leftrightarrow$ \PL$[\bmu]$, and $\mu$ is pseudo-tight if and only if $\bmu$ is tight.
\item $\bmu$ is continuous from below if and only if $\mu$ is continuous from above; $\bmu$ is continuous from above if and only if $\mu$ is continuous from below.
\item $\bmu$ is order continuous if and only if $\mu$ satisfies \CBX.
\item $\bmu$ is strongly order continuous if and only if $B_n\upto B$ and $\mu(B)=\mu(X)$ imply $\mu(B_n)\to\mu(X)$.
\item $\bmu$ is exhaustive if and only if $\mu(X\setminus E_n)\to\mu(X)$ for every sequence $(E_n)$ of pairwise disjoint sets.
\item $\bmu$ is weakly null-additive if and only if $\mu$ satisfies
\DW: $\mu(G)=\mu(H)=\mu(X)$ implies $\mu(G\cap H)=\mu(X)$.
\item $\bmu$ has the p.g.p.\ if and only if $\mu$ satisfies \DP: for every $\varepsilon>0$ there is $\delta>0$ such that $\mu(G)>\mu(X)-\delta$ and $\mu(H)>\mu(X)-\delta$ imply $\mu(G\cap H)>\mu(X)-\varepsilon$.
\item $\bmu$ has property (S) if and only if $\mu$ satisfies \DS: whenever $\mu(B_n)\to\mu(X)$, some subsequence satisfies $\mu\bigl(\bigcup_k\bigcap_{i\ge k}B_{n_i}\bigr)=\mu(X)$.
\item $\bmu$ fulfils condition (E) if and only if $\mu$ satisfies \DE: for every double sequence with $B^{(m)}_n\upto B^{(m)}$ as $n\to\infty$ for each $m$ and $\mu\bigl(\bigcap_m B^{(m)}\bigr)=\mu(X)$, there are increasing sequences $(n_i)$, $(m_i)$ with $\lim_{k}\mu\bigl(\bigcap_{i\ge k}B^{(m_i)}_{n_i}\bigr)=\mu(X)$.
\item $\bmu$ is regular if and only if $\mu$ satisfies \PR: for every $A\in\cB$ and $\varepsilon>0$ there are $F\in\cC$, $G\in\cO$ with $F\subseteq A\subseteq G$ and $\mu(F\cup G^c)>\mu(X)-\varepsilon$.
\end{enumerate}
\end{proposition}

\begin{proof}
The map $A\mapsto A^c$ is a bijection of $\cB$ onto itself which reverses inclusions, turns non-decreasing sequences into non-increasing ones, and interchanges unions and intersections. Moreover, $\nu(A)=0$ if and only if $\mu(A^c)=\mu(X)$, and $\nu(A)<\eta$ if and only if $\mu(A^c)>\mu(X)-\eta$.\\
(a) The identity is the definition of $\nu$. Hence "$\mu(F)>\mu(X)-\varepsilon$" and "$\nu(X\setminus F)<\varepsilon$" are equivalent for the same set $F$, which gives the first two equivalences and the statement on tightness; the third equivalence follows by applying the first one to $\bmu$ and using $\overline{\bmu}=\mu$.\\
(b) If $A_n\upto A$, then $A_n^c\downto A^c$ and $\nu(A_n)=\mu(X)-\mu(A_n^c)$, so $\nu(A_n)\to\nu(A)$ if and only if $\mu(A_n^c)\to\mu(A^c)$. As $(A_n)$ runs through all non-decreasing sequences, $(A_n^c)$ runs through all non-increasing ones; since $\mu$ is finite, the side condition in continuity from above is automatic. The second equivalence is the first one applied to $\bmu$.\\
(c) $A_n\downto\emptyset$ if and only if $B_n:=A_n^c\upto X$, and $\nu(A_n)\to0$ if and only if $\mu(B_n)\to\mu(X)$.\\
(d) $A_n\downto A$ with $\nu(A)=0$ if and only if $B_n:=A_n^c\upto B:=A^c$ with $\mu(B)=\mu(X)$.\\
(e) $\nu(E_n)=\mu(X)-\mu(X\setminus E_n)$.\\
(f), (g) Put $G=E^c$ and $H=F^c$; then $E\cup F$ corresponds to $G\cap H$.\\
(h) $\nu\bigl(\bigcap_k\bigcup_{i\ge k}A_{n_i}\bigr)=0$ if and only if $\mu\bigl(\bigcup_k\bigcap_{i\ge k}A_{n_i}^c\bigr)=\mu(X)$.\\
(i) Take complements in the definition of condition (E) for $\nu$.\\
(j) $(G\setminus F)^c=F\cup G^c$.
\end{proof}

\begin{remark}[Comparison with known "pseudo" conditions]\label{rem:pseudo-comparison}
Li et. al. \cite{LOM21} introduced absolute continuity of types E, R and L for a pair $(\lambda,\nu)$ of monotone measures, written $\lambda\ll_E\nu$, $\lambda\ll_R\nu$, $\lambda\ll_L\nu$; for $\lambda=\nu$ these reduce to condition (E), property (S) and strong order continuity of $\nu$, respectively \cite[Proposition~3.4]{LOM21}. In this language:
\begin{itemize}[leftmargin=1em]
\item The pseudo-condition (E) of \cite[Definition~4.6]{LMPK15}, denoted $[\overline{E}]$ in \cite[Definition~2.1 and Note~2.4]{LOM21}, has the hypothesis $\mu\bigl(\bigcup_m E^{(m)}\bigr)=0$ and the same conclusion as \DE{} written for $E=B^c$; it is $\bmu\ll_E\mu$, whereas \DE{} is condition (E) for $\bmu$.
\item Similarly, the pseudo-property (S), denoted $[\overline{S}]$ in \cite[Definition~5.19]{LMPK15} and \cite[Definition~2.2]{LOM21}, is $\bmu\ll_R\mu$, whereas \DS{} is property (S) for $\bmu$.
\item Strong order pseudo-continuity \cite[Definition~2.3]{LOM21} is $\bmu\ll_L\mu$, whereas \cref{prop:dictionary}(d) is strong order continuity of $\bmu$, that is, $\bmu\ll_L\bmu$.
\end{itemize}
Thus the ``pseudo'' conditions in the literature are mixed conditions ({} \cite[Proposition~3.4]{LOM21}), while \DE, \DS{} and \cref{prop:dictionary}(d) are the corresponding pure conditions for $\bmu$.
\end{remark}

\begin{proposition}[The condition (CS)]\label{prop:CS}
Let $v$ be a finite monotone measure.
\begin{enumerate}[(a)]
\item $v$ satisfies \CS{} if and only if $\bv$ is a submeasure.
\item Every supermodular $v$ satisfies \CS. In particular, belief functions and $2$-monotone capacities satisfy \CS.
\item Let $\cP$ be a nonempty family of finitely additive set functions $p\colon\cB\to[0,\infty)$ with a common total mass $p(X)=c$. Then $\inf_{p\in\cP}p$ satisfies \CS. In particular, every coherent lower probability satisfies \CS.
\item If $v$ satisfies \CS, then $v\le\bv$.
\item \CS{} is strictly weaker than supermodularity.
\item Every finite submeasure $\mu$ satisfies $\bmu\le\mu$.
\end{enumerate}
\end{proposition}

\begin{proof}
(a) The inequality $\bv(A\cup B)\le\bv(A)+\bv(B)$ is equivalent to $v(A^c)+v(B^c)\le v(X)+v(A^c\cap B^c)$, and $(A,B)\mapsto(A^c,B^c)$ is a bijection of $\cB\times\cB$.\\
(b) $v(A)+v(B)\le v(A\cup B)+v(A\cap B)\le v(X)+v(A\cap B)$.\\
(c) For each $p\in\cP$ we have $v(A)+v(B)\le p(A)+p(B)=p(A\cup B)+p(A\cap B)\le c+p(A\cap B)$, where $v=\inf_{p\in\cP}p$; now take the infimum over $p$.\\
(d) Take $B=A^c$ in \CS.\\
(e) This is the familiar phenomenon of coherent lower probabilities which are not $2$-monotone; here is a concrete instance. Let $X=\{a,b,c,d\}$ with the discrete metric, $p=(2/3,1/3,0,0)$, $q=(1/3,0,1/3,1/3)$ and $v=\min(p,q)$, which satisfies \CS{} by (c). For $A=\{a,c\}$ and $B=\{b,c\}$ we have $v(A)=2/3$, $v(B)=1/3$, $v(A\cup B)=\min(1,2/3)=2/3$ and $v(A\cap B)=0$, so $v(A\cup B)+v(A\cap B)=2/3<1=v(A)+v(B)$.\\
(f) $\mu(X)\le\mu(A)+\mu(A^c)$.
\end{proof}
\noindent
The following facts are well known; we include the short proof for convenience.

\begin{lemma}\label{lem:submeasure}
Let $\mu$ be a finite submeasure. The following are equivalent: (i) $\mu$ is order continuous; (ii) $\mu$ is continuous from above; (iii) $\mu$ is continuous. If they hold, $\mu$ is countably subadditive. Moreover, every finite submeasure is null-additive (hence weakly null-additive), has the p.g.p.\ (with $\delta=\varepsilon/2$), and is autocontinuous from above and from below.
\end{lemma}

\begin{proof}
(iii)$\Rightarrow$(ii)$\Rightarrow$(i) is trivial, since $\mu$ is finite. For (i)$\Rightarrow$(iii), let first $A_n\downto A$; then $\mu(A)\le\mu(A_n)\le\mu(A)+\mu(A_n\setminus A)$ and $A_n\setminus A\downto\emptyset$. If $A_n\upto A$, then $\mu(A_n)\le\mu(A)\le\mu(A_n)+\mu(A\setminus A_n)$ and $A\setminus A_n\downto\emptyset$. Countable subadditivity follows by applying finite subadditivity to $\bigcup_{i\le n}A_i$ and letting $n\to\infty$, using continuity from below. The remaining claims follow from $\mu(E)\le\mu(E\cup F)\le\mu(E)+\mu(F)$ and from $\mu(E)\ge\mu(E\setminus F)\ge\mu(E)-\mu(F)$, the latter because $E\subseteq(E\setminus F)\cup F$.
\end{proof}

\begin{corollary}\label{cor:CS-continuity}
Let $v$ be a finite monotone measure satisfying \CS. Then \CBX, continuity from below and continuity of $v$ are equivalent.
\end{corollary}

\begin{proof}
By \cref{prop:CS}(a), $\bv$ is a finite submeasure. By \cref{prop:dictionary}(b),(c), \CBX{} for $v$ is order continuity of $\bv$, continuity from below of $v$ is continuity from above of $\bv$, and continuity of $v$ is continuity of $\bv$. Now apply \cref{lem:submeasure} to $\bv$.
\end{proof}

The condition \CS{} cannot be dropped in \cref{cor:CS-continuity}: the set function $w$ of \cref{ex:vacuous} and the set function $\nu$ of \cref{ex:meagre} both satisfy \CBX{} and are continuous from below, but they are not continuous from above.

\section{Necessary conditions for Lusin's theorem}\label{sec:necessity}

\begin{theorem}\label{thm:necessity}
Let $\mu$ be a monotone measure on a metric space $X$.
\begin{enumerate}[(a)]
\item If \Lu$[\mu]$ holds, then $\mu$ is regular.
\item If \LK$[\mu]$ holds, then $\mu$ is order continuous and tight.
\item If $X$ is compact and \Lu$[\mu]$ holds, then $\mu$ is regular, order continuous and exhaustive.
\item Let $\mu$ be finite. If \PL$[\mu]$ holds, then $\mu$ satisfies \PR. If \PLK$[\mu]$ holds, then $\mu$ satisfies \CBX{} and is pseudo-tight. If $X$ is compact and \PL$[\mu]$ holds, then $\mu$ satisfies \PR{} and \CBX.
\end{enumerate}
\end{theorem}

\begin{proof}
(a) Fix $A\in\cB$ and $\varepsilon>0$, and put $f:=\1_A$. By \Lu$[\mu]$ there is $F\in\cC$ with $\mu(X\setminus F)<\varepsilon$ such that $f|_F$ is continuous. The sets $A\cap F=(f|_F)^{-1}([1/2,\infty))$ and $F\setminus A=(f|_F)^{-1}((-\infty,1/2])$ are closed in $F$, hence closed in $X$. Put $F':=A\cap F\subseteq A$ and $G:=X\setminus(F\setminus A)$; then $F'$ is closed, $G$ is open and $G\supseteq A$. We claim that $G\setminus F'=X\setminus F$. Indeed, a point $x\notin F$ belongs neither to $F\setminus A$ nor to $A\cap F$, so $x\in G\setminus F'$; a point $x\in F$ lies either in $A\cap F=F'$ or in $F\setminus A$, i.e.\ $x\notin G$. Hence $\mu(G\setminus F')=\mu(X\setminus F)<\varepsilon$.\\
(b) Tightness follows by applying \LK$[\mu]$ to $f\equiv0$. For order continuity, let $A_n\downto\emptyset$ and put $f:=\sum_n\1_{A_n}$. Since the sets $A_n$ decrease, $x\in A_n$ if and only if $n<N(x):=\min\{n: x\notin A_n\}$, and this minimum exists because $\bigcap_nA_n=\emptyset$. Thus $f(x)=N(x)-1$ is finite, and $\{f\ge k\}=A_k$ for $k\ge1$, so $f$ is Borel. Given $\varepsilon>0$, choose $K\in\cK$ with $\mu(X\setminus K)<\varepsilon$ such that $f|_K$ is continuous. The sets $C_k:=A_k\cap K=(f|_K)^{-1}([k,\infty))$ are closed in $K$, hence compact; they decrease and $\bigcap_kC_k=\emptyset$. By the finite intersection property, $C_N=\emptyset$ for some $N$, that is, $A_N\subseteq X\setminus K$. Hence $\mu(A_n)\le\mu(A_N)\le\mu(X\setminus K)<\varepsilon$ for all $n\ge N$.\\
(c) On a compact space \Lu{} and \LK{} coincide, so (a) and (b) apply. For exhaustivity, let $(E_n)$ be pairwise disjoint; then $\bigcup_{k\ge n}E_k\downto\emptyset$, and $\mu(E_n)\le\mu\bigl(\bigcup_{k\ge n}E_k\bigr)\to0$.\\
(d) Apply (a)--(c) to $\bmu$ and use \cref{prop:dictionary}(a),(c),(j).
\end{proof}

\begin{remark}\label{rem:necessity}
(i) Compact Lusin sets are needed in \cref{thm:necessity}(b); see \cref{ex:infinite-sets}.\\
(ii) Part (a) is elementary and possibly folklore. However, we have not found it stated in \cite{Cao15,JS96,JWZ99,LMPK15,LM11,LY04,WYT13}; in these works regularity is used as a tool for proving the Lusin property rather than obtained as a consequence of it.\\
(iii) The compactness argument in part (b) is the one used by Wiesel \cite[Proposition~5.4]{Wie22}, where it is written for sub-additive capacities (normalized and continuous from below) and decreasing sequences of open sets. Part (b) is the corresponding statement for arbitrary monotone measures and arbitrary decreasing sequences of Borel sets. For uncertain measures on compact metric spaces, related necessity results were obtained by Yang et. al. \cite{YZH24}.
\end{remark}

\section{Lusin's theorem for submeasures}\label{sec:submeasures}

We first recall that an order continuous finite submeasure has the Lusin property. This is a special case of several known theorems.

\begin{theorem}\label{thm:sub-known}
Let $\mu$ be a finite order continuous submeasure on a metric space $X$. Then \Lu$[\mu]$ holds.
\end{theorem}

\begin{proof}
By \cref{lem:submeasure}, $\mu$ is finite and continuous, hence continuous from below, exhaustive (by order continuity) and autocontinuous from above and from below; it is also weakly null-additive and has the p.g.p.. Therefore, each of the following results applies to $\mu$:
\begin{itemize}[leftmargin=1em]
\item \cite[Theorem~5]{JS96} (a fuzzy measure, i.e.\ monotone and continuous from above and from below, which is $\sigma$-finite, exhaustive and autocontinuous from above; the conclusion is stated with $\mu(X\setminus F)\le\varepsilon$, which suffices). In \cite[Corollary~4]{JWZ99} this result is quoted with the word ``autocontinuous''; both readings are covered by \cref{lem:submeasure}.
\item \cite[Theorem~5]{JWZ99} (a lower semicontinuous fuzzy measure, i.e.\ monotone and continuous from below, which is exhaustive and has the p.g.p.).
\item \cite[Theorem~4]{LY04} (a finite continuous weakly null-additive fuzzy measure).
\item \cite[Corollary~5.2]{LM11} (continuous from below, exhaustive, p.g.p.).
\item \cite[Theorem~5.3]{LM11} (condition (E) and the p.g.p.); condition (E) holds by \cite[Proposition~3.1]{LM11}, since $\mu$ is finite and continuous.
\end{itemize}
A self-contained proof, which follows the same route as \cite[Theorem~1]{JS96} and \cite[Proposition~5.2 and Theorem~5.3]{Wie22}, is given in \cref{app:A1,app:A2}.
\end{proof}

\begin{corollary}\label{cor:compact-sub}
Let $X$ be a compact metric space and $\mu$ a finite submeasure on $X$. The following are equivalent:
\begin{enumerate}[(i)]
\item \Lu$[\mu]$;
\item $\mu$ is order continuous;
\item $\mu$ is continuous from above;
\item $\mu$ is continuous.
\end{enumerate}
\end{corollary}

\begin{proof}
(i)$\Rightarrow$(ii) is \cref{thm:necessity}(c), (ii)$\Rightarrow$(i) is \cref{thm:sub-known}, and (ii)$\Leftrightarrow$(iii)$\Leftrightarrow$(iv) is \cref{lem:submeasure}.
\end{proof}

\begin{remark}\label{rem:Wiesel}
Wiesel \cite[Theorem~1.2]{Wie22} works on a compact metric space with normalized sub-additive \emph{capacities}, which by his definition are continuous from below. His characterizing condition is continuity from above along decreasing sequences of \emph{open} sets:
\[
\text{(W)}\quad O_n\in\cO,\ O_n\downto O\ \Longrightarrow\ v(O_n)\to v(O).
\]
Under continuity from below, (W) is equivalent to continuity from above \cite[Remark~4.1]{Wie22}. The two results are not comparable: \cref{cor:compact-sub} does not contain \cite[Theorem~1.2]{Wie22}, because without continuity from below we have not been able to show that (W) suffices; conversely, \cite[Theorem~1.2]{Wie22} does not contain \cref{cor:compact-sub}, because it assumes continuity from below. What \cref{cor:compact-sub} adds is that \Lu{} forces full continuity, so that Wiesel's standing assumption of continuity from below holds automatically for any finite submeasure with \Lu; and that, without this assumption, \Lu{} is equivalent to order continuity along all Borel sequences. Whether (W) alone suffices for finite submeasures on compact spaces is left open (\cref{sec:concl}, Problem~\ref{Q2}).
\end{remark}

\begin{theorem}\label{thm:compact-Lusin}
Let $\mu$ be a finite submeasure on a metric space $X$. Then we have
\begin{enumerate}[(a)]
\item \LK$[\mu]$ holds if and only if $\mu$ is order continuous and tight.
\item If $X$ is Polish, every order continuous finite submeasure is tight. Hence \LK$[\mu]$ $\Leftrightarrow$ $\mu$ is order continuous $\Leftrightarrow$ $\mu$ is continuous.
\end{enumerate}
\end{theorem}

\begin{proof}
(a) The implication "$\Rightarrow$" is \cref{thm:necessity}(b). Conversely, let $f$ be Borel and $\varepsilon>0$. By \cref{thm:sub-known} there is $F\in\cC$ with $\mu(X\setminus F)<\varepsilon/2$ and $f|_F$ continuous, and by tightness there is $K\in\cK$ with $\mu(X\setminus K)<\varepsilon/2$. Then $F\cap K\in\cK$, $f|_{F\cap K}$ is continuous, and $\mu(X\setminus(F\cap K))\le\mu(X\setminus F)+\mu(X\setminus K)<\varepsilon$.\\
(b) The implication "$\Leftarrow$" is known: \cite[Theorem~6]{JWZ99} applies, since $\mu$ is continuous from below, exhaustive and has the p.g.p.\ (\cref{lem:submeasure}); so do \cite[Theorem~6]{JS96} (finite, exhaustive, autocontinuous from above) and \cite[Theorem~3.3]{Cao15} (finite, continuous, weakly null-additive). The tightness statement is contained in \cite[Theorem~2]{JWZ99}, \cite[Theorem~2]{JS96}, \cite[Lemma~3.2]{Cao15} and \cite[Theorem~17]{WYT13}; a self-contained argument in the spirit of Ulam's theorem is given in \cref{app:A3}. The implication ``$\Rightarrow$'' is \cref{thm:necessity}(b), and the last equivalence is \cref{lem:submeasure}.
\end{proof}

The following standard example separates \Lu{} from order continuity and from \LK, even on a Polish space.

\begin{example}\label{ex:infinite-sets}
Let $X=\N$ with the discrete metric (a Polish space), and let $\mu(A)=1$ if $A$ is infinite and $\mu(A)=0$ otherwise. Then we have
\begin{itemize}[leftmargin=1em]
\item $\mu$ is a finite submeasure, since $A\cup B$ infinite implies that $A$ or $B$ is infinite;
\item $\mu$ is not order continuous: $\{n,n+1,\dots\}\downto\emptyset$, but each of these sets has measure $1$;
\item \Lu$[\mu]$ holds with $F=X$, because every function on a discrete space is continuous;
\item \LK$[\mu]$ fails: compact subsets of $\N$ are finite, so $\mu(X\setminus K)=1$ for every $K\in\cK$.
\end{itemize}
\end{example}

\section{Pseudo-Lusin theorems}\label{sec:pseudo}

\subsection{The pseudo-Lusin theorem under (CS)}

\begin{theorem}\label{thm:CS}
Let $v$ be a finite monotone measure satisfying \CS{}. Then we have
\begin{enumerate}[(a)]
\item Let $X$ be any metric space. If $v$ satisfies \CBX{} (equivalently, $v$ is continuous from below, or continuous; see \cref{cor:CS-continuity}), then \PL$[v]$ and \Lu$[v]$ hold with the same closed set $F$:
\[
v(F)>v(X)-\varepsilon\;\;\text{and}\;\; v(X\setminus F)\le\bv(X\setminus F)<\varepsilon .
\]
\item If $X$ is compact, then \PL$[v]$ $\Leftrightarrow$ \CBX{} $\Leftrightarrow$ $v$ is continuous from below $\Leftrightarrow$ $v$ is continuous.
\item Let $X$ be any metric space. Then \PLK$[v]$ holds if and only if $v$ satisfies \CBX{} and is pseudo-tight.
\item If $X$ is Polish, then \PLK$[v]$ $\Leftrightarrow$ \CBX{} $\Leftrightarrow$ $v$ is continuous from below $\Leftrightarrow$ $v$ is continuous.
\end{enumerate}
\end{theorem}

\begin{proof}
By \cref{prop:CS}(a), $\bv$ is a finite submeasure, and by \cref{prop:dictionary}(c), \CBX{} for $v$ is order continuity of $\bv$.\\
(a) Apply \cref{thm:sub-known} to $\bv$, and then use \cref{prop:dictionary}(a) and \cref{prop:CS}(d).\\
(b) If \PL$[v]$ holds, then \Lu$[\bv]$ holds by \cref{prop:dictionary}(a), so $\bv$ is order continuous by \cref{thm:necessity}(c). The converse is (a), and the remaining equivalences are \cref{cor:CS-continuity}.\\
(c) Apply \cref{thm:compact-Lusin}(a) to $\bv$ and translate by \cref{prop:dictionary}(a),(c).\\
(d) Apply \cref{thm:compact-Lusin}(b) to $\bv$, and then use \cref{cor:CS-continuity}.
\end{proof}

\begin{remark}\label{rem:CS-status}
Part (a) is a consequence of \cref{thm:sub-known} by duality. In part (d) the implication "$\Leftarrow$" is known in dual form, namely \cite[Theorem~6]{JWZ99} applied to $\bv$; the necessity parts in (b)--(d) appear to be new. No continuity from above of $v$ is assumed in \cref{thm:CS}; on compact or Polish spaces it is a consequence of \PL{} or \PLK, respectively.
\end{remark}

\subsection{Comparison of (L) and (PL)}

\begin{proposition}\label{prop:L-vs-PL}
Let $\mu$ be a finite monotone measure. Then we have
\begin{enumerate}[(a)]
\item If $\mu\le\bmu$ (in particular, if $\mu$ satisfies \CS), then \PL$[\mu]\Rightarrow$ \Lu$[\mu]$ and \PLK$[\mu]\Rightarrow$ \LK$[\mu]$.
\item If $\bmu\le\mu$ (in particular, if $\mu$ is a submeasure), then \Lu$[\mu]\Rightarrow$ \PL$[\mu]$ and \LK$[\mu]\Rightarrow$ \PLK$[\mu]$.
\item Neither converse holds (\cref{ex:vacuous}).
\end{enumerate}
\end{proposition}

\begin{proof}
(a) If $\mu(F)>\mu(X)-\varepsilon$, then by \cref{prop:dictionary}(a) we get $\mu(X\setminus F)\le\bmu(X\setminus F)=\mu(X)-\mu(F)<\varepsilon$.\\
(b) is symmetric: $\mu(X\setminus F)<\varepsilon$ gives $\mu(X)-\mu(F)=\bmu(X\setminus F)\le\mu(X\setminus F)<\varepsilon$. The parenthetical cases follow from \cref{prop:CS}(d),(f).
\end{proof}

\subsection{Pseudo-Lusin theorems obtained from known Lusin theorems}\label{subsec:routine}

For a finite monotone measure $\mu$, \PL$[\mu]$ follows from any Lusin theorem applied to $\bmu$. Translating the hypotheses by \cref{prop:dictionary}, one obtains the following statements; they are routine, but we record them for later reference.

\begin{proposition}\label{prop:PL-list}
Let $\mu$ be a finite monotone measure on a metric space $X$. Then each of the conditions \textup{(PL1)}--\textup{(PL5)} in \cref{tab:PL} implies \PL$[\mu]$.
\begin{table}[H]
\centering
\small
\caption{Sufficient conditions for \PL$[\mu]$, $\mu$ finite. The third column lists the Lusin theorem applied to $\bmu$, and the last column the hypotheses of that theorem for $\bmu$, obtained via \cref{prop:dictionary}.}\label{tab:PL}
\begin{tabular}{@{}l>{\raggedright\arraybackslash}p{0.33\textwidth}>{\raggedright\arraybackslash}p{0.22\textwidth}>{\raggedright\arraybackslash}p{0.27\textwidth}@{}}
\toprule
 & Hypotheses on $\mu$ & Source & Hypotheses for $\bmu$ \\
\midrule
(PL1) & $\mu$ continuous and \DW & \cite[Theorem~4]{LY04} & finite, continuous (b), w.n.a.\ (f) \\
(PL2) & $\mu$ continuous from above; $\mu(X\setminus E_n)\to\mu(X)$ for pairwise disjoint $(E_n)$; \DP & \cite[Corollary~5.2]{LM11}, \cite[Theorem~5]{JWZ99} & continuous from below (b), exhaustive (e), p.g.p.\ (g) \\
(PL3) & \DE{} and \DP & \cite[Theorem~5.3]{LM11} & (E) (i), p.g.p.\ (g) \\
(PL4) & [$B_n\upto B$, $\mu(B)=\mu(X)$ $\Rightarrow$ $\mu(B_n)\to\mu(X)$], \DS{} and \DW & \cite[Corollary~5.3]{LM11} & strongly order continuous (d), (S) (h), w.n.a.\ (f) \\
(PL5) & $\mu$ continuous from below, \DE{} and \DW & \cref{app:A4}; also \cite[Theorem~6.7]{LMPK15} & continuous from above (b), (E) (i), w.n.a.\ (f) \\
\bottomrule
\end{tabular}
\end{table}
\end{proposition}

\begin{proof}
In each row of \cref{tab:PL}, the hypotheses on $\mu$ are equivalent, by the indicated parts of \cref{prop:dictionary}, to the hypotheses on $\bmu$ in the last column, and $\bmu$ is finite. The cited theorem gives \Lu$[\bmu]$, which is \PL$[\mu]$ by \cref{prop:dictionary}(a). For (PL5), the finite case of \cite[Theorem~6.7]{LMPK15} (which is (L5) in \cite[Section~7]{LM11}) is proved in \cref{app:A4}.
\end{proof}

Compact-set versions on Polish spaces follow in the same way from \cite[Theorem~6]{JWZ99}, \cite[Theorem~6]{JS96} and \cite[Theorem~3.3]{Cao15}. Note that supermodularity gives \DW{} directly: if $v(A)=v(B)=v(X)$, then $v(A\cap B)\ge v(A)+v(B)-v(A\cup B)\ge v(X)$. Moreover, \CS{} implies both \DW{} and \DP, because $\bv$ is then a submeasure (\cref{prop:CS}(a), \cref{lem:submeasure}).

\subsection{Examples}

\begin{example}[The vacuous lower probability on {$[0,1]$}]\label{ex:vacuous}
Let $X=[0,1]$, let $v(A)=1$ if $A=X$ and $v(A)=0$ otherwise, and let $w:=\bv$, so that $w(A)=1$ if and only if $A\ne\emptyset$.\\
\emph{Properties of $v$:}
\begin{enumerate}[1.]
\item $v$ is monotone, $v(\emptyset)=0$ and $v(X)=1$.
\item $v$ is supermodular. If $A=X$, then $v(A\cup B)+v(A\cap B)=1+v(B)=v(A)+v(B)$; the case $B=X$ is symmetric, and otherwise the right-hand side of the supermodularity inequality is $0$. In fact $v$ is the belief function with basic probability assignment $m(X)=1$; in particular \CS{} holds.
\item $v$ is continuous from above, hence order continuous and strongly order continuous. Indeed, let $A_n\downto A$. If $A=X$, then $A_n=X$ for all $n$. Otherwise pick $x\notin A$; then $x\notin A_N$ for some $N$, and $v(A_n)=0$ for $n\ge N$.
\item $v$ is neither continuous from below nor \CBX: $A_n=\{0\}\cup[1/n,1]\upto X$, but $v(A_n)=0$ because $1/(2n)\notin A_n$.
\item $v$ is not weakly null-additive, hence does not have the p.g.p.: $v([0,1/2])=v((1/2,1])=0$ but $v(X)=1$.
\item $v$ fulfils condition (E). Let $E^{(m)}_n\downto E^{(m)}$ for each $m$, with $D:=\bigcup_mE^{(m)}$ and $v(D)=0$, i.e.\ $D\ne X$; pick $x_0\notin D$. For each $m$ we have $x_0\notin E^{(m)}_n$ for all large $n$. Choose $m_i=i$ and $n_1<n_2<\cdots$ with $x_0\notin E^{(i)}_{n_i}$. Then $\bigcup_{i\ge k}E^{(m_i)}_{n_i}\ne X$, so its $v$-value is $0$.
\item $v$ does not have property (S): let $(q_n)$ enumerate $\Q\cap[0,1]$ and $A_n=X\setminus\{q_n\}$. Then $v(A_n)=0$, but $\bigcap_k\bigcup_{i\ge k}A_{n_i}=X$ for every subsequence.
\item $v$ is regular: for $A\notin\{\emptyset,X\}$ pick $x\notin A$ and take $F=\emptyset$, $G=X\setminus\{x\}$; then $v(G\setminus F)=0$. (For $A=\emptyset$ or $A=X$ take $F=G=A$.)
\item \Lu$[v]$ and \LK$[v]$ hold, but for a trivial reason: $F=\{0\}$ is compact and $v(X\setminus F)=v((0,1])=0$. This is as degenerate as \cite[Example~2]{LY04} and \cite[Example~7.1]{LM11}, but on a non-discrete space.
\item \PL$[v]$ fails: take $f=\1_{[0,1/2]}$ and $0<\varepsilon\le1$. Then $v(F)>1-\varepsilon\ge0$ forces $F=X$, and $f$ is not continuous on $X$.
\end{enumerate}
\noindent
\emph{Properties of $w$:} $w$ is a submeasure and is continuous from below. It is not order continuous, since $(0,1/n)\downto\emptyset$ while $w((0,1/n))=1$. It is weakly null-additive and has the p.g.p.\ (\cref{lem:submeasure}). \PL$[w]$ holds, since it is equivalent to \Lu$[v]$ (\cref{prop:dictionary}(a)), while \Lu$[w]$ fails by \cref{thm:necessity}(c). Finally, $w$ is not regular: for $A=[0,1/2]$, the inequality $w(G\setminus F)<1$ with $F\subseteq A\subseteq G$, $F$ closed and $G$ open, forces $G=F=A$, so $[0,1/2]$ would be clopen in the connected space $[0,1]$.\\
\emph{Consequences:} \Lu{} does not imply \PL, even for totally monotone capacities; \PL{} does not imply \Lu, even for submeasures which are continuous from below; and, for measures satisfying \CS, continuity from above alone does not give \PL.
\end{example}

As far as we know, the following example has not been used before in connection with Lusin's theorem.

\begin{example}[Meagre and comeagre sets]\label{ex:meagre}
Let $X=[0,1]$ (compact and Polish). Put $\mu(A)=1$ if $A$ is comeagre and $\mu(A)=0$ otherwise, and let $\nu:=\bmu$. Then $\nu(A)=1-\mu(X\setminus A)=1$ if and only if $X\setminus A$ is not comeagre, that is, if and only if $A$ is non-meagre; thus $\nu(A)=0$ exactly when $A$ is meagre.\\
We use two facts throughout. (Baire) By the Baire category theorem, every nonempty open subset of $[0,1]$ is non-meagre; in particular $[0,1]$ is not meagre. Moreover, meagre sets form a $\sigma$-ideal. (N) A closed set with empty interior is nowhere dense, hence meagre; therefore, a closed non-meagre set has nonempty interior.\\
\emph{Properties of $\mu$:}
\begin{enumerate}[1.]
\item \emph{$\mu$ is a normalized monotone measure.} $\mu(\emptyset)=0$, since otherwise $X$ would be meagre, contrary to (Baire); $\mu(X)=1$; and if $A\subseteq B$ with $A$ comeagre, then $X\setminus B\subseteq X\setminus A$ is meagre.
\item \emph{$\mu$ is supermodular, hence satisfies \CS.} If $A$ and $B$ are both comeagre, then so is $A\cap B$ (its complement is the union of two meagre sets), and the left-hand side $\mu(A\cup B)+\mu(A\cap B)$ equals $2$. If exactly one of them, say $A$, is comeagre, then $A\cup B$ is comeagre and the left-hand side is $1+\mu(A\cap B)\ge1=\mu(A)+\mu(B)$. Otherwise the right-hand side is $0$.
\item \emph{$\mu$ is continuous from above} along all Borel sequences. Let $A_n\downto A$. If every $A_n$ is comeagre, then $X\setminus A=\bigcup_n(X\setminus A_n)$ is meagre, so $\mu(A)=1=\mu(A_n)$. Otherwise some $A_N$ is not comeagre; then no $A_n$ with $n\ge N$ is comeagre, and neither is $A\subseteq A_N$, so $\mu(A_n)=0=\mu(A)$ for $n\ge N$.
\item \emph{$\mu$ is order continuous, strongly order continuous and exhaustive.} The first two properties follow from 3. For exhaustivity, note that two disjoint sets cannot both be comeagre (the union of their complements, which is $X$, would be meagre); hence, for pairwise disjoint $(E_n)$, $\mu(E_n)=0$ for all but at most one $n$.
\item \emph{$\mu$ is neither continuous from below nor \CBX.} Indeed $A_n=\{0\}\cup[1/n,1]\upto X$, but $X\setminus A_n=(0,1/n)$ is a nonempty open set, hence non-meagre by (Baire), so $\mu(A_n)=0$.
\item \emph{$\mu$ is not regular; hence \Lu$[\mu]$ fails} by \cref{thm:necessity}(a). Take $A=\Q\cap[0,1]$. A closed set $F\subseteq A$ is countable, so its interior is empty (nonempty open subsets of $[0,1]$ are uncountable), and $F$ is meagre by (N). An open set $G\supseteq A$ is dense, so $X\setminus G$ is closed with empty interior, hence meagre. Therefore, $X\setminus(G\setminus F)=(X\setminus G)\cup F$ is meagre, $G\setminus F$ is comeagre, and $\mu(G\setminus F)=1$.\\
The failure of \Lu$[\mu]$ can also be seen directly: take $f=\1_{\Q}$ and $0<\varepsilon\le1$. If $F$ is closed and $\mu(X\setminus F)<1$, then $X\setminus F$ is not comeagre, so $F$ is not meagre and, by (N), contains a nondegenerate interval $I$. But $f|_I$ is discontinuous at every point of $I$, since both $\Q$ and its complement are dense.
\item \PL$[\mu]$ fails as well: by \cref{prop:L-vs-PL}(a), \PL{} implies \Lu{} under \CS; alternatively, apply \cref{thm:CS}(b), since $\mu$ does not satisfy \CBX.
\item \emph{$\mu$ is not a submeasure}: $\mu([0,1/2])=\mu((1/2,1])=0<1=\mu(X)$; neither half is comeagre, since each complement contains a nonempty open interval.
\end{enumerate}
\noindent 
\emph{Properties of $\nu=\bmu$:}
\begin{enumerate}[1.]\setcounter{enumi}{8}
\item \emph{$\nu$ is a submeasure}: if $A\cup B$ is non-meagre, then $A$ or $B$ is non-meagre, since the union of two meagre sets is meagre.
\item \emph{$\nu$ is continuous from below, hence \CBX.} Let $A_n\upto A$. If $A$ is meagre, then so is every $A_n$. If $A$ is non-meagre, then some $A_n$ is non-meagre, since a countable union of meagre sets is meagre; thus $\nu(A_n)=1=\nu(A)$ for all large $n$.
\item \emph{$\nu$ is neither order continuous nor continuous from above}: $(0,1/n)\downto\emptyset$, while $\nu((0,1/n))=1$.
\item \emph{\PL$[\nu]$ fails}, since \PL$[\nu]$ is equivalent to \Lu$[\mu]$ by \cref{prop:dictionary}(a), and \Lu$[\mu]$ fails by 6. \emph{\Lu$[\nu]$ fails as well}, by \cref{thm:necessity}(c), because $\nu$ is not order continuous and $X$ is compact.
\item \emph{$\nu$ does not satisfy \CS}: by \cref{prop:CS}(a), \CS{} for $\nu$ means that $\overline{\nu}=\mu$ is a submeasure, which fails by 8.
\end{enumerate}
\noindent 
\emph{Consequences:}
\begin{itemize}[leftmargin=2em]
\item \CS{} cannot be omitted in \cref{thm:CS}(a),(b), nor in \cref{cor:CS-continuity}: $\nu$ satisfies \CBX{} and is continuous from below on a compact space, yet \PL$[\nu]$ fails.
\item Subadditivity is essential in \cref{cor:compact-sub}: $\mu$ is finite, order continuous, continuous from above, exhaustive and even supermodular on a compact space, yet \Lu$[\mu]$ fails.
\item Compare \cref{ex:vacuous}, where $v$ is also supermodular and continuous from above, does not satisfy \CBX, and nevertheless has \Lu{} (trivially). Hence, for supermodular capacities, continuity from above neither forces nor excludes \Lu.
\item The measure $w$ of \cref{ex:vacuous} satisfies \CBX{} but not \CS{} and has \PL, whereas $\nu$ satisfies \CBX{} but not \CS{} and fails \PL. Hence \CBX{} alone neither forces nor excludes \PL.
\item $\mu$ is \emph{not} a Dempster capacity in the sense of \cite{CM02} (see \cref{sec:CM}): condition (iv) there fails for $A_n=[1/n,1]$, since $A_n\subseteq\inter A_{n+1}=(1/(n+1),1]$ and $\bigcup_nA_n=(0,1]$ is comeagre, while $\mu(A_n)=0$. So this example gives no new information on \cite[Theorem~1]{CM02}.
\item The example does not answer Problem~\ref{Q1} of \cref{sec:concl}, since $\mu$ is not continuous from below, nor Problem~\ref{Q3}, since $\mu$ is not regular.
\end{itemize}
\end{example}

\section{On the theorem of Castaldo and Marinacci}\label{sec:CM}

Castaldo and Marinacci \cite[Definition~1]{CM02} call a set function $v\colon\cB\to[0,1]$ on the Borel sets of a topological space a \emph{Dempster capacity} if (i) $v(\emptyset)=0$ and $v(X)=1$; (ii) $v$ is monotone; (iii) $v(A_n)\to v\bigl(\bigcap_nA_n\bigr)$ for every non-increasing sequence $(A_n)$; and (iv) $v(A_n)\to v\bigl(\bigcup_nA_n\bigr)$ for every non-decreasing sequence with $A_n\subseteq\inter A_{n+1}$. A \emph{Choquet capacity} in \cite{CM02} satisfies (i)--(iii) and continuity from below along increasing sequences of open sets. The main result \cite[Theorem~1]{CM02} asserts: \emph{if $v$ is a supermodular Dempster capacity on a Polish space $X$ and $u\colon X\to\R$ is Borel, then for every $\varepsilon>0$ there is a compact set $K_\varepsilon$ with $v(K_\varepsilon)\ge1-\varepsilon$ such that $u|_{K_\varepsilon}$ is continuous.}

\subsection{A counterexample and the corrected statement}

\begin{example}\label{ex:CM-counter}
Let $X=[0,1]$ and let $v$ be the vacuous capacity of \cref{ex:vacuous}. Then $v$ satisfies all conditions of \cite[Definition~1]{CM02}: (i) and (ii) are clear; (iii) is \cref{ex:vacuous}, item 3; for (iv), let $(A_n)$ be non-decreasing with $A_n\subseteq\inter A_{n+1}$ and $A:=\bigcup_nA_n$. If $A\ne X$, then $v(A_n)=0=v(A)$. If $A=X$, then $(\inter A_{n+1})_n$ is an increasing open cover of the compact space $X$, so $X=\inter A_{N+1}\subseteq A_{N+1}$ for some $N$, and $v(A_n)=1$ for $n\ge N+1$. Moreover, $v$ is supermodular, even totally monotone (\cref{ex:vacuous}, item 2). However, for $u=\1_{[0,1/2]}$ and $\varepsilon\in(0,1)$, a compact set $K$ with $v(K)\ge1-\varepsilon$ must be $K=X$, and $u$ is not continuous on $X$. Thus \cite[Theorem~1]{CM02} fails for $v$. The obstruction is precisely the one identified in \cref{thm:CM-corrected} below: $v$ is not continuous from below (\cref{ex:vacuous}, item 4).
\end{example}

\begin{remark}\label{rem:Wiesel-credit}
Wiesel \cite[Section~3]{Wie22} was the first to refute \cite[Theorem~1]{CM02}, whose statement he recalls in dual form. He used a Huber--Strassen type capacity \cite{HS73}, $v_W(A)=\min\bigl(Q(A^{\delta})+\varepsilon_W,1\bigr)$ for nonempty closed $A$, where $Q$ is a Borel probability measure and $A^{\delta}$ denotes the closed $\delta$-neighbourhood of $A$, together with the function $u=\1_{\{0\}}$; any $\varepsilon_W,\delta>0$ are allowed. For $\varepsilon_W=1$ one gets $v_W(A)=1$ for every nonempty closed $A$, which agrees on closed sets with $w=\bv$. Hence \cref{ex:CM-counter} is, in dual form, the degenerate case $\varepsilon_W=1$ of Wiesel's example. It adds only two small points: the axioms of \cite{CM02} are verified on all Borel sets, directly in the original supermodular formulation and without appeal to \cite{HS73}; and, as observed in \cref{ex:random-sets} below, $v(A)=P(\Gamma\subseteq A)$ is the distribution of the constant compact-valued random set $\Gamma\equiv[0,1]$, so \cite[Theorem~1]{CM02} fails already within the class of compact-valued random correspondences which motivated it.
\end{remark}

\begin{theorem}[Corrected Castaldo--Marinacci theorem]\label{thm:CM-corrected}
Let $X$ be a Polish space and $v\colon\cB\to[0,1]$ a monotone set function with $v(\emptyset)=0$ and $v(X)=1$ which is supermodular (or merely satisfies \CS). The following are equivalent:
\begin{enumerate}[(i)]
\item for every Borel function $u\colon X\to\R$ and every $\varepsilon>0$ there is a compact set $K_\varepsilon$ with $v(K_\varepsilon)\ge1-\varepsilon$ such that $u|_{K_\varepsilon}$ is continuous (the conclusion of \cite[Theorem~1]{CM02});
\item $v$ is continuous from below along all non-decreasing sequences of Borel sets or, equivalently, along those increasing to $X$;
\item $v$ is continuous from above and from below.
\end{enumerate}
\end{theorem}

\begin{proof}
Condition (i) is \PLK$[v]$ (\cref{rem:Lusin-basic}(ii)), so the equivalences follow from \cref{thm:CS}(d).
\end{proof}

\begin{remark}\label{rem:CM-meaning}
\Cref{thm:CM-corrected} shows that \cite[Theorem~1]{CM02} becomes correct if condition (iv) of \cite[Definition~1]{CM02} (continuity from below along sequences with $A_n\subseteq\inter A_{n+1}$) is replaced by continuity from below along all non-decreasing Borel sequences; condition (iii) of \cite[Definition~1]{CM02} is then automatic, and on Polish spaces the strengthened hypothesis is also necessary. The sufficiency part is, in dual form, \cite[Theorem~6]{JWZ99} applied to $w:=\bv$; what is new is the necessity and the ``if and only if'' formulation.

In Wiesel's dual language, let $w:=\bv$, a normalized sub-additive set function. On a Polish space, $w$ has the compact-set Lusin property if and only if $w$ is continuous from above along \emph{all} non-increasing Borel sequences, and no continuity from below of $w$ is assumed. We do not regard this as a strict extension of \cite[Theorem~1.2]{Wie22}: the sufficiency on Polish spaces predates \cite{Wie22} (it is in \cite{JWZ99}); Wiesel himself points to the non-compact case in \cite[Lemma~4.3]{Wie22}; and his condition (W) concerns open sets only, the two formulations being equivalent only under continuity from below ({} \cref{rem:Wiesel}).
\end{remark}

\subsection{On Lemma 3 of Castaldo--Marinacci and the gap in the proof}\label{subsec:gap}

In this subsection we make Wiesel's diagnosis \cite[Section~3]{Wie22} of the gap in the proof of \cite[Theorem~1]{CM02} more explicit. No novelty is claimed for the diagnosis itself. For a set function $v$ with $v(X)=1$, let $\core(v)$ be the set of countably additive Borel probability measures $p$ with $p\ge v$.

The proof in \cite{CM02} relies on \cite[Lemma~3]{CM02}, which, for a metric space $X$ and a \emph{regular capacity} $v$, meaning that
\begin{equation}\label{eq:CM-regular}
v(A)=\sup\{v(C): C\subseteq A,\ C\in\cC\}=\inf\{v(G): G\supseteq A,\ G\in\cO\}\quad(A\in\cB),
\end{equation}
asserts that $\core(v)$ is weak$^*$ compact, and that, if in addition $X$ is Polish, for each $\varepsilon>0$ there is a compact $K_\varepsilon$ with $p(K_\varepsilon)\ge1-\varepsilon$ for all $p\in\core(v)$. We call \eqref{eq:CM-regular} \emph{CM-regularity}. The word ``capacity'' alone is not defined in \cite{CM02}; its Definition~1 and the subsequent paragraph define only Dempster and Choquet capacities. Two readings are therefore possible, and the status of the lemma depends on the reading.

\begin{proposition}[Literal reading]\label{prop:R1}
If ``capacity'' means a normalized monotone set function, both parts of \cite[Lemma~3]{CM02} are false. Indeed, let $X=\N$ with the discrete metric (a Polish space) and $v(A)=1$ if $A=\N$, $v(A)=0$ otherwise. Then $v$ is supermodular and CM-regular, but $\core(v)$ is neither tight nor compact.
\end{proposition}

\begin{proof}
Every subset of $\N$ is clopen, so \eqref{eq:CM-regular} holds with $C=G=A$; supermodularity is checked exactly as in \cref{ex:vacuous}, item 2. The core $\core(v)$ is the set of all probability measures on $\N$; it contains every Dirac measure $\delta_n$. Every compact $K\subseteq\N$ is finite, so $\delta_n(K)=0$ for $n\notin K$, and $\core(v)$ is not tight.
Nor is $\core(v)$ compact, in any reasonable reading of ``weak$^*$''. For the topology of weak convergence of probability measures: if a subsequence $(\delta_{n_j})$ converged weakly to some $p$, then for each finite $F\subseteq\N$, which is clopen with empty boundary, $p(F)=\lim_j\delta_{n_j}(F)=0$, whence $p(\N)=0$, a contradiction. (Alternatively, on a Polish space relative weak compactness implies tightness, by the converse part of Prokhorov's theorem \cite{Bil68}.) For the topology $\sigma(\mathrm{ca},B_b)$ induced by bounded Borel functions: if a subnet $(\delta_{n_\alpha})$ converged to a countably additive $p$, then $p(F)=0$ for every finite $F$, hence $p(\N)=0$ by countable additivity, whereas $p(\N)=\lim_\alpha\delta_{n_\alpha}(\N)=1$.
\end{proof}

\begin{remark}
In the proof of the first part of \cite[Lemma~3]{CM02}, a limit obtained from the Alaoglu theorem need not be countably additive; the argument tacitly assumes that it is. The second part is derived in \cite{CM02} from the first by a classical result in \cite{Bil68}. Note also that the set function $v$ of \cref{prop:R1} is not a Choquet capacity: the open sets $\{1,\dots,n\}$ increase to $\N$, but have $v$-value $0$.
\end{remark}

\begin{proposition}\label{prop:R2}
Let $X$ be a Polish space and $v\colon\cB\to[0,1]$ monotone with $v(X)=1$ and continuous from below along increasing sequences of open sets (for instance, a Choquet capacity in the sense of \cite{CM02}). Then $\core(v)$ is tight. If, in addition, $v$ is CM-regular, then $\core(v)$ is compact in the topology of weak convergence of probability measures.
\end{proposition}

\begin{proof}
Let $\rho$ be a complete metric compatible with the topology of $X$, and let $(x_i)$ be dense in $X$. Fix $\varepsilon>0$. For $k\ge1$ the open sets $U_{k,N}:=\bigcup_{i\le N}B_\rho(x_i,1/k)$ increase to $X$ as $N\to\infty$, so by continuity from below along open sets we may choose $N_k$ with $v(U_{k,N_k})>1-\varepsilon2^{-k}$. Let $K:=\bigcap_k\cl U_{k,N_k}$. Then $K$ is closed, hence complete, and for each $k$,
\(
K\subseteq\cl U_{k,N_k}\subseteq\bigcup_{i\le N_k}B_\rho(x_i,2/k),
\)
so $K$ is totally bounded and therefore compact. For $p\in\core(v)$, countable subadditivity of $p$ and $p\ge v$ give
\[
p(X\setminus K)\le\sum_kp(X\setminus\cl U_{k,N_k})\le\sum_k\bigl(1-p(U_{k,N_k})\bigr)\le\sum_k\bigl(1-v(U_{k,N_k})\bigr)<\varepsilon .
\]
Hence $\core(v)$ is tight. \\
By Prokhorov's theorem \cite{Bil68}, $\core(v)$ is relatively compact in the space $P(X)$ of Borel probability measures with the topology of weak convergence, which is metrizable. It remains to show that $\core(v)$ is closed. Let $p_j\in\core(v)$ converge weakly to $p$. For closed $C$, the portmanteau theorem gives $p(C)\ge\limsup_jp_j(C)\ge v(C)$. For $A\in\cB$, inner CM-regularity gives $v(A)=\sup\{v(C):C\subseteq A,\,C\in\cC\}\le\sup\{p(C):C\subseteq A,\,C\in\cC\}\le p(A)$. Hence $p\in\core(v)$, and $\core(v)$ is compact.
\end{proof}
\noindent
If ``weak$^*$'' in \cite{CM02} refers to a different topology, the compactness statement of \cref{prop:R2} should be read with the topology of weak convergence of $P(X)$. Note that tightness does not use regularity, whereas compactness does.

In either reading, the proof of \cite[Theorem~1]{CM02} breaks down when \cite[Lemma~3]{CM02} is applied to a finer topology $\tau_u$. Under the literal reading the lemma itself is false (\cref{prop:R1}). Under the Choquet-capacity reading its hypothesis fails in $\tau_u$, as the following proposition shows for the capacity of \cref{ex:CM-counter}.

\begin{proposition}\label{prop:tau-u}
Let $X=[0,1]$ with its usual topology $\tau$, let $v$ be as in \cref{ex:CM-counter}, and let $u=\1_{[0,1/2]}$. Let $\tau_u$ be any Polish topology on $[0,1]$ which is finer than $\tau$, has the same Borel sets and makes $u$ continuous (such a topology is provided by \cite[Lemma~1]{CM02}; {} \cite{Kec95,Sri98}). Then $v$ is not a Choquet capacity with respect to $\tau_u$, and $\core(v)$ is not $\tau_u$-tight.
\end{proposition}

\begin{proof}
The sets $[0,1/2]=u^{-1}(\{1\})$ and $(1/2,1]=u^{-1}(\{0\})$ are $\tau_u$-closed, hence $\tau_u$-clopen. The sets $A_n=[0,1/2]\cup(1/2+1/n,1]$ are $\tau_u$-open (since $\tau\subseteq\tau_u$) and increase to $X$, but $v(A_n)=0$; so $v$ is not a Choquet capacity with respect to $\tau_u$.\\
The core $\core(v)$ consists of all Borel probability measures, and in particular contains the Dirac measures $\delta_{1/2+1/n}$. Let $K$ be $\tau_u$-compact. Since the identity $(X,\tau_u)\to(X,\tau)$ is continuous, $K$ is $\tau$-compact; moreover, $K\cap(1/2,1]$ is a $\tau_u$-closed subset of $K$, hence $\tau_u$-compact, hence $\tau$-compact. Thus $K\cap(1/2,1]$ is a compact subset of $(1/2,1]$, so it is bounded away from $1/2$, and $\delta_{1/2+1/n}(K)=0$ for all large $n$. Hence no $\tau_u$-compact set has $p$-measure $\ge1-\varepsilon$ for all $p\in\core(v)$ when $\varepsilon<1$.
\end{proof}

This is precisely Wiesel's observation in \cite[Section~3]{Wie22} that $v$ need not be a capacity with respect to $\tau_u$, so that the arguments of \cite[Lemma~3]{CM02} do not apply and tightness of the core fails. Our additions are only the observation that the lemma, read literally, is false while it holds on Polish spaces under the Choquet-capacity reading, and the verification for every admissible $\tau_u$.

\begin{remark}[Regularity]\label{rem:CM-regularity}
Let $v$ be a finite monotone measure satisfying \CS. If $\bv$ is regular (in the sense of \cref{def:regular}), then $v$ is CM-regular. Indeed, let $A\in\cB$ and $F\subseteq A\subseteq G$ with $F\in\cC$ and $G\in\cO$. Applying \CS{} to $A$ and $X\setminus(G\setminus F)$, and using $A\setminus(G\setminus F)=F$, we obtain $v(A)-v(F)\le\bv(G\setminus F)$. Applying \CS{} to $G$ and $X\setminus(G\setminus A)$, and using $G\setminus(G\setminus A)=A$, we obtain $v(G)-v(A)\le\bv(G\setminus A)\le\bv(G\setminus F)$.
\end{remark}

\subsection{A natural positive class: random finite sets}

\begin{example}\label{ex:random-sets}
Let $(\Omega,\mathscr{A},P)$ be a probability space, let $N\colon\Omega\to\{1,2,\dots\}$ be $\mathscr{A}$-measurable, and let $\xi_1,\xi_2,\allowbreak\ldots\colon\Omega\to X$ be $\mathscr{A}/\cB$-measurable. Put $\Gamma(\omega):=\{\xi_1(\omega),\allowbreak\dots,\xi_{N(\omega)}(\omega)\}$ and $v(A):=P(\Gamma\subseteq A)$ for $A\in\cB$. Then
\begin{itemize}[leftmargin=1em]
\item $v$ is well defined, since $\{\Gamma\subseteq A\}=\bigcup_n\bigl(\{N=n\}\cap\bigcap_{i\le n}\{\xi_i\in A\}\bigr)\in\mathscr{A}$.
\item $v(\emptyset)=0$, $v(X)=1$, and $v$ is monotone.
\item $v$ is supermodular, because $\{\Gamma\subseteq A\}\cap\{\Gamma\subseteq B\}=\{\Gamma\subseteq A\cap B\}$ and $\{\Gamma\subseteq A\}\cup\{\Gamma\subseteq B\}\subseteq\{\Gamma\subseteq A\cup B\}$, so that $v(A\cup B)\ge v(A)+v(B)-v(A\cap B)$.
\item $v$ is continuous from below: if $A_n\upto A$ and $\Gamma(\omega)\subseteq A$, then the finitely many points of $\Gamma(\omega)$ lie in $A_n$ for all large $n$.
\end{itemize}
By \cref{thm:CS}, \PL$[v]$ holds on every metric space $X$, and \PLK$[v]$ holds on every Polish space. For $N\equiv1$ this is the classical Lusin theorem. In contrast, the constant compact random set $\Gamma\equiv[0,1]$ gives the capacity of \cref{ex:CM-counter}; and the constant random set $\Gamma\equiv\Q\cap[0,1]$ gives $v(A_n)=0$ for $A_n=X\setminus\{q_k:k>n\}\upto X$, where $(q_k)$ enumerates $\Q\cap[0,1]$, so that \CBX{} fails and \PL$[v]$ fails on $[0,1]$ by \cref{thm:CS}(b).
\end{example}

\section{Transfer between two monotone measures}\label{sec:pairs}

Since the Lusin property involves a single set function, statements for pairs of monotone measures, in the spirit of \cite{LOM21}, reduce to transfer results. These are elementary. We use the notation $\lambda\ll_I\nu$, $\lambda\ll_{VI}\nu$, $\lambda\ll_E\nu$, $\lambda\ll_R\nu$, $\lambda\ll_L\nu$ of \cite[Definitions~3.1--3.5]{LOM21}; in particular, $\lambda\ll_{VI}\nu$ means that $\nu(A_n)\to0$ implies $\lambda(A_n)\to0$, and $\lambda\ll_E\nu$ means that for every double sequence with $A^{(m)}_n\downto A^{(m)}$ as $n\to\infty$ and $\nu\bigl(\bigcup_mA^{(m)}\bigr)=0$ there are increasing $(n_i)$, $(m_i)$ with $\lim_k\lambda\bigl(\bigcup_{i\ge k}A^{(m_i)}_{n_i}\bigr)=0$.

\begin{proposition}\label{prop:transfer}
Let $\lambda,\nu$ be monotone measures on $X$ such that $\nu(O_n)\to0$ implies $\lambda(O_n)\to0$ for every sequence $(O_n)$ of open sets. Then \Lu$[\nu]\Rightarrow$ \Lu$[\lambda]$ and \LK$[\nu]\Rightarrow$ \LK$[\lambda]$. For finite $\lambda,\nu$, the analogous statement holds for \PL{} and \PLK, with $\bla$, $\bnu$ in place of $\lambda$, $\nu$ in the hypothesis.
\end{proposition}

\begin{proof}
First, for every $\varepsilon>0$ there is $\delta>0$ such that $O\in\cO$ and $\nu(O)<\delta$ imply $\lambda(O)<\varepsilon$; otherwise there would be $\varepsilon_0>0$ and open sets $O_n$ with $\nu(O_n)<1/n$ and $\lambda(O_n)\ge\varepsilon_0$, contrary to the hypothesis. Now, given a Borel function $f$ and $\varepsilon>0$, take $F$ (closed or compact) with $\nu(X\setminus F)<\delta$ and $f|_F$ continuous; since $X\setminus F$ is open, $\lambda(X\setminus F)<\varepsilon$. The statement for \PL{} and \PLK{} follows by applying the first part to $\bla$ and $\bnu$ and using \cref{prop:dictionary}(a).
\end{proof}

\begin{lemma}\label{lem:S-E}
If $\nu$ has property (S) and $\lambda\ll_E\nu$, then $\lambda\ll_{VI}\nu$. Consequently, \Lu$[\nu]\Rightarrow$ \Lu$[\lambda]$.
\end{lemma}

\begin{proof}
Suppose that $\nu(A_n)\to0$ but $\lambda(A_n)\ge\varepsilon_0>0$ for all $n$. By property (S) there is a subsequence $(B_i)$ of $(A_n)$ such that $\nu(D)=0$, where $D:=\bigcap_j\bigcup_{i\ge j}B_i$. For every $m$ put $C^{(m)}_j:=\bigcup_{i\ge j}B_i$; then $C^{(m)}_j\downto D$ as $j\to\infty$, and $\nu\bigl(\bigcup_mD\bigr)=\nu(D)=0$. Since $\lambda\ll_E\nu$, there are increasing sequences $(j_r)$ and $(m_r)$ with $\lambda\bigl(\bigcup_{r\ge k}C^{(m_r)}_{j_r}\bigr)\to0$. This union equals $\bigcup_{i\ge j_k}B_i\supseteq B_{j_k}$, hence $\lambda(B_{j_k})\to0$, a contradiction. The last assertion follows from \cref{prop:transfer}.
\end{proof}

\begin{example}\label{ex:pairs}
Let $\nu=m$ be the Lebesgue measure on $X=[0,1]$ and $\lambda(A)=1$ if $m(A)>0$, $\lambda(A)=0$ otherwise. Then $\lambda\ll_Im$. Also $\lambda\ll_Rm$: if $m(A_n)\to0$, choose a subsequence with $\sum_im(A_{n_i})<\infty$ and apply the Borel--Cantelli lemma. On the other hand, $\lambda$ is neither $\ll_Em$, nor $\ll_Lm$, nor $\ll_{VI}m$, as the sequence $(0,1/n)$ shows. Finally, \Lu$[m]$ holds but \Lu$[\lambda]$ fails: for $f=\1_{[0,1/2]}$ and $\varepsilon<1$, a closed $F$ with $\lambda(X\setminus F)<\varepsilon$ satisfies $m(X\setminus F)=0$; then $X\setminus F$ is an open null set, hence empty, so $F=X$, and $f$ is not continuous on $X$.
\end{example}

\section{An application: approximation by continuous functions}\label{sec:approx}

The following is a routine analogue of \cite[Theorems~5 and~6]{LY04}, \cite[Theorem~5.4]{LM11} and \cite[Theorems~4.1 and~4.2]{Cao15}. For a monotone measure $\mu$ and a nonnegative Borel function $g$ we use the Sugeno integral $(Su)\!\int g\,d\mu=\sup_{t\ge0}\min\bigl(t,\mu(\{g\ge t\})\bigr)$ and the Choquet integral $(C)\!\int g\,d\mu=\int_0^\infty\mu(\{g>t\})\,dt$. Following \cite{LMPK15}, a sequence $(f_n)$ converges to $f$ \emph{pseudo-in measure} $\mu$ if $\mu(\{|f_n-f|<\alpha\})\to\mu(X)$ for every $\alpha>0$.

\begin{theorem}\label{thm:approx}
Let $X$ be a metric space, $v$ a finite monotone measure satisfying \CS{} and \CBX, and $f$ a Borel function. Then there are continuous functions $\varphi_n$ on $X$ with the following properties. Then
\begin{enumerate}[(a)]
\item $\bv(\{|\varphi_n-f|\ge\alpha\})<1/n$ for every $\alpha>0$. In particular $v(\{|\varphi_n-f|<\alpha\})\to v(X)$, i.e.\ $(\varphi_n)$ converges to $f$ pseudo-in measure $v$ \textup{(} \cite[Proposition~3.5]{LMPK15}\textup{)}; moreover $v(\{|\varphi_n-f|\ge\alpha\})\to0$.
\item $(Su)\!\int|\varphi_n-f|\,dv\le(Su)\!\int|\varphi_n-f|\,d\bv\le1/n$.
\item If $|f|\le M$, one may take $|\varphi_n|\le M$, and then $(C)\!\int|\varphi_n-f|\,dv\le(C)\!\int|\varphi_n-f|\,d\bv\le2M/n$.
\end{enumerate}
\end{theorem}

\begin{proof}
By \cref{thm:CS}(a), for each $n$ there is $F_n\in\cC$ with $\bv(X\setminus F_n)<1/n$ and $f|_{F_n}$ continuous. By Tietze's extension theorem for closed subsets of metric spaces, $f|_{F_n}$ extends to a continuous function $\varphi_n$ on $X$; if $|f|\le M$, we replace $\varphi_n$ by $(\varphi_n\wedge M)\vee(-M)$. For $t>0$ we have $\{|\varphi_n-f|\ge t\}\subseteq X\setminus F_n$, which gives (a), since $v(\{|\varphi_n-f|<\alpha\})=v(X)-\bv(\{|\varphi_n-f|\ge\alpha\})$. For the Sugeno integral, $(Su)\!\int g\,d\bv=\sup_t\min\bigl(t,\bv(\{g\ge t\})\bigr)\le1/n$ with $g=|\varphi_n-f|$. For the Choquet integral in (c), the integrand $\bv(\{g>t\})$ is less than $1/n$ for $0<t<2M$ and vanishes for $t\ge2M$. The comparisons with $v$ follow from $v\le\bv$ (\cref{prop:CS}(d)).
\end{proof}

\section{Concluding remarks and open problems}\label{sec:concl}

We have studied Lusin's theorem for monotone measures from the point of view of necessary conditions and of conjugate duality. The necessity results of \cref{thm:necessity} turn several known sufficient conditions for submeasures into characterizations, and the pseudo-Lusin property \PL{} provides a convenient framework for set functions arising in imprecise probability, for which large lower probability of the Lusin set is the natural requirement. In particular, the Lusin-type theorem of Castaldo and Marinacci holds, on Polish spaces, exactly for those supermodular capacities which are continuous from below. We hope that these observations may be of some use to researchers working on regularity and convergence theorems for non-additive measures.

Several questions remain open.
\begin{enumerate}[label=\textbf{Q\arabic*.},ref=Q\arabic*,leftmargin=3em]
\item\label{Q1} (Li and Yasuda \cite[Section~7]{LY04}.) Does every finite continuous monotone measure on a compact metric space have the Lusin property? By duality, this is equivalent to asking whether every finite continuous monotone measure has \PL. \Cref{ex:meagre} does not answer this question, since its $\mu$ is not continuous from below. We note that any monotone measure $\lambda$ with $\nu(O_n)\to0\Rightarrow\lambda(O_n)\to0$ for open sets, where $\nu$ has \Lu, has \Lu{} as well (\cref{prop:transfer}).
\item\label{Q2} Let $X$ be compact and $\mu$ a finite submeasure satisfying Wiesel's open-set condition (W), but not assumed to be continuous from below. Does \Lu$[\mu]$ hold? Equivalently, does (W) imply order continuity for such $\mu$?
\item\label{Q3} On a compact metric space, does ``regular and order continuous'' imply \Lu{} for an arbitrary finite monotone measure? (The measure $\mu$ of \cref{ex:meagre} is not regular.)
\item\label{Q4} (Li and Mesiar \cite[Section~7]{LM11}.) Are condition (E) and weak null-additivity sufficient for \Lu?
\item\label{Q5} Do \Lu$[\lambda]$ and \Lu$[\nu]$ imply \Lu$[\lambda\vee\nu]$? This is true on compact spaces for finite submeasures, by \cref{cor:compact-sub}, since $\lambda\vee\nu$ is then again an order continuous submeasure.
\item\label{Q6} Does $\lambda\ll_E\nu$ alone allow one to transfer \Lu{} from $\nu$ to $\lambda$?
\end{enumerate}

\appendix
\section{Self-contained proofs}\label{app:proofs}

In \cref{app:A1,app:A2,app:A3} we give, for completeness, proofs of known results for submeasures; \cref{app:A4} proves the row (PL5) of \cref{tab:PL}. Throughout \cref{app:A1,app:A2,app:A3}, $\mu$ is a finite order continuous submeasure on a metric space $X$; by \cref{lem:submeasure} it is continuous and countably subadditive.

\subsection{Regularity}\label{app:A1}
(\cite[Theorem~1]{JS96}, \cite[Theorem~1]{JWZ99}, \cite[Proposition~5.2]{Wie22}.) Let $\cM$ be the family of all $A\in\cB$ such that for every $\varepsilon>0$ there are $F\in\cC$ and $G\in\cO$ with $F\subseteq A\subseteq G$ and $\mu(G\setminus F)<\varepsilon$. We show that $\cM=\cB$. Indeed,
\begin{itemize}[leftmargin=1em]
\item \emph{Closed sets belong to $\cM$.} For $F\in\cC$, the sets $G_m:=\{x:\dist(x,F)<1/m\}$ are open and $G_m\setminus F\downto\emptyset$, so $\mu(G_m\setminus F)\to0$ by order continuity.
\item \emph{$\cM$ is closed under complements.} If $F\subseteq A\subseteq G$, then $G^c\subseteq A^c\subseteq F^c$ and $F^c\setminus G^c=G\setminus F$.
\item \emph{$\cM$ is closed under countable unions.} Let $A_n\in\cM$ and choose $F_n\subseteq A_n\subseteq G_n$ with $\mu(G_n\setminus F_n)<\varepsilon2^{-n-1}$. Put $G=\bigcup_nG_n$, $\Phi=\bigcup_nF_n$ and $\Phi_N=\bigcup_{n\le N}F_n$. Since $G\setminus\Phi\subseteq\bigcup_n(G_n\setminus F_n)$, countable subadditivity gives $\mu(G\setminus\Phi)\le\varepsilon/2$. As $G\setminus\Phi_N\downto G\setminus\Phi$, continuity from above yields $\mu(G\setminus\Phi_N)<\varepsilon$ for some $N$, and $\Phi_N$ is closed.
\end{itemize}
Hence $\cM$ is a $\sigma$-algebra containing the closed sets, so $\cM=\cB$.

\subsection{Lusin's theorem with closed sets}\label{app:A2}
First let $f=\sum_na_n\1_{A_n}$, where $(A_n)$ is a countable Borel partition of $X$. By \cref{app:A1}, choose closed $F_n\subseteq A_n$ with $\mu(A_n\setminus F_n)<\varepsilon2^{-n-2}$, and, by order continuity, choose $N$ with $\mu\bigl(X\setminus\bigcup_{n\le N}A_n\bigr)<\varepsilon/2$. Put $F=\bigcup_{n\le N}F_n$. By finite subadditivity, $\mu(X\setminus F)<3\varepsilon/4$. The sets $F_n$, $n\le N$, are finitely many pairwise disjoint closed sets, so each of them is relatively open in $F$; since $f$ is constant on each $F_n$, $f|_F$ is locally constant, hence continuous.\\
For general $f$, let $g_k=2^{-k}\lfloor2^kf\rfloor$, which is countably valued. By the first step take closed sets $H_k$ with $\mu(X\setminus H_k)<\varepsilon2^{-k-1}$ and $g_k|_{H_k}$ continuous, and put $F=\bigcap_kH_k$. By countable subadditivity, $\mu(X\setminus F)\le\varepsilon/2$. Since $g_k\to f$ uniformly on $F$, $f|_F$ is continuous.

\subsection{Tightness on Polish spaces}\label{app:A3}
Let $X$ be Polish, let $\rho$ be a complete metric compatible with the topology of $X$ (so that open, closed and compact sets are the same for $\rho$ and $d$), and let $(x_i)$ be dense in $X$. Put $U_{k,N}=\bigcup_{i\le N}\overline{B}_\rho(x_i,1/k)$. Since $U_{k,N}\upto X$ as $N\to\infty$, order continuity allows us to choose $N_k$ with $\mu(X\setminus U_{k,N_k})<\varepsilon2^{-k}$. The set $K=\bigcap_kU_{k,N_k}$ is closed and totally bounded (for each $k$ it is covered by $N_k$ balls of radius $2/k$), hence compact, and $\mu(X\setminus K)<\varepsilon$ by countable subadditivity.

\subsection{Lusin's theorem under continuity from above, condition (E) and weak null-additivity (finite case)}\label{app:A4}

\begin{proposition}\label{prop:A4}
Let $\lambda$ be a finite monotone measure on a metric space $X$ which is continuous from above, fulfils condition (E) and is weakly null-additive. Then \Lu$[\lambda]$ holds.
\end{proposition}
This is the finite case of \cite[Theorem~6.7]{LMPK15}, which is (L5) in \cite[Section~7]{LM11}. The only external ingredients are the following two results of Li and Mesiar:
\begin{itemize}[leftmargin=2em]
\item \cite[Proposition~3.2]{LM11}: condition (E) implies order continuity;
\item \cite[Proposition~3.7, (1)$\Rightarrow$(2)]{LM11}: weak null-additivity and condition (E) imply
\begin{itemize}
\item[($\star$)] for every $\varepsilon>0$ and every double sequence $(E^{(m)}_n)$ with $E^{(m)}_n\downto E^{(m)}$ as $n\to\infty$ and $\lambda(E^{(m)})=0$ for each $m$, there are indices $n_m$ with $\lambda\bigl(\bigcup_mE^{(m)}_{n_m}\bigr)<\varepsilon$.
\end{itemize}
\end{itemize}
Finiteness of $\lambda$ is used only through continuity from above, whose side condition is then automatic.
\begin{proof}
\emph{Step 1 (regularity).} Let $\cM$ be defined as in \cref{app:A1}, with $\lambda$ in place of $\mu$.
\begin{itemize}[leftmargin=2em]
\item \emph{Closed sets belong to $\cM$.} For $F\in\cC$, with $G_m$ as in \cref{app:A1}, $G_m\setminus F\downto\emptyset$; since $\lambda$ is order continuous by \cite[Proposition~3.2]{LM11}, $\lambda(G_m\setminus F)\to0$.
\item \emph{$\cM$ is closed under complements}, exactly as in \cref{app:A1}.
\item \emph{$\cM$ is closed under countable unions.} Let $A_n\in\cM$. For each $n$ and $i$ choose $F_{n,i}\in\cC$ and $G_{n,i}\in\cO$ with $F_{n,i}\subseteq A_n\subseteq G_{n,i}$ and $\lambda(G_{n,i}\setminus F_{n,i})<1/i$, and put $F^{(k)}_n:=\bigcup_{i\le k}F_{n,i}$ and $G^{(k)}_n:=\bigcap_{i\le k}G_{n,i}$. Since $G^{(k)}_n\setminus F^{(k)}_n\subseteq G_{n,k}\setminus F_{n,k}$, this replacement only shrinks the differences. Thus $F^{(k)}_n$ is closed, $G^{(k)}_n$ is open, $F^{(k)}_n\subseteq A_n\subseteq G^{(k)}_n$, $F^{(k)}_n\upto_k$, $G^{(k)}_n\downto_k$ and $\lambda(G^{(k)}_n\setminus F^{(k)}_n)<1/k$. The sets $E^{(n)}_k:=G^{(k)}_n\setminus F^{(k)}_n$ decrease in $k$ to $D_n:=\bigcap_kE^{(n)}_k$, and $\lambda(D_n)\le1/k$ for all $k$, so $\lambda(D_n)=0$. Given $\varepsilon>0$, ($\star$) yields indices $k_n$ with $\lambda\bigl(\bigcup_nE^{(n)}_{k_n}\bigr)<\varepsilon$. Put $G:=\bigcup_nG^{(k_n)}_n$ (open, containing $\bigcup_nA_n$), $\Phi:=\bigcup_nF^{(k_n)}_n$ and $\Phi_N:=\bigcup_{n\le N}F^{(k_n)}_n$ (closed, contained in $\bigcup_nA_n$). Every $x\in G\setminus\Phi$ lies in some $G^{(k_n)}_n$ but not in $F^{(k_n)}_n$, so $G\setminus\Phi\subseteq\bigcup_nE^{(n)}_{k_n}$ and $\lambda(G\setminus\Phi)<\varepsilon$. Since $G\setminus\Phi_N\downto G\setminus\Phi$, continuity from above gives $\lambda(G\setminus\Phi_N)<\varepsilon$ for some $N$.
\end{itemize}
Hence $\cM=\cB$, that is, $\lambda$ is regular.\\
\emph{Step 2 (countably valued functions).} Let $g=\sum_na_n\1_{A_n}$, where $(A_n)$ is a countable Borel partition of $X$. By Step 1, and by taking finite unions, choose closed sets $F^{(k)}_n\subseteq A_n$ with $F^{(k)}_n\upto_k$ and $\lambda(A_n\setminus F^{(k)}_n)<1/k$; then $A_n\setminus F^{(k)}_n$ decreases in $k$ to a $\lambda$-null set. Given $\delta>0$, ($\star$) gives indices $k_n$ with $\lambda\bigl(\bigcup_n(A_n\setminus F^{(k_n)}_n)\bigr)<\delta$. Since $(A_n)$ is a partition, $X\setminus\bigcup_nF^{(k_n)}_n=\bigcup_n(A_n\setminus F^{(k_n)}_n)$. As $X\setminus\bigcup_{n\le N}F^{(k_n)}_n\downto X\setminus\bigcup_nF^{(k_n)}_n$, continuity from above gives $N$ with $\lambda(X\setminus H)<\delta$, where $H:=\bigcup_{n\le N}F^{(k_n)}_n$.\\
Call a set \emph{admissible for $g$} if it is of the form $\bigcup_{n\le N}C_n$ with closed $C_n\subseteq A_n$. An admissible set is closed, and $g$ is continuous on it, because the $C_n$ are finitely many pairwise disjoint closed sets, so each of them is relatively open, and $g$ is constant on each $C_n$. A finite union of admissible sets is again admissible (group the pieces according to $n$). We have shown: for every $\delta>0$ there is an admissible set $H$ with $\lambda(X\setminus H)<\delta$.\\
\emph{Step 3 (general $f$).} Let $g_k:=2^{-k}\lfloor2^kf\rfloor$; these are countably valued and $0\le f-g_k<2^{-k}$. For each $k$ and $j$, Step 2 gives a set admissible for $g_k$ whose complement has $\lambda$-measure less than $1/j$; let $H^{(j)}_k$ be the union of the first $j$ of these sets. Then $H^{(j)}_k$ is admissible for $g_k$, increasing in $j$, and $\lambda(X\setminus H^{(j)}_k)<1/j$. Hence $X\setminus H^{(j)}_k\downto_j D_k$ with $\lambda(D_k)=0$. Given $\varepsilon>0$, ($\star$) gives indices $j_k$ with $\lambda\bigl(\bigcup_k(X\setminus H^{(j_k)}_k)\bigr)<\varepsilon$. The set $F:=\bigcap_kH^{(j_k)}_k$ is closed and $\lambda(X\setminus F)<\varepsilon$. Each $g_k$ is continuous on $F\subseteq H^{(j_k)}_k$, and $g_k\to f$ uniformly, so $f|_F$ is continuous.
\end{proof}
\noindent
\emph{Application to (PL5).} Let $\mu$ be finite and continuous from below, and let it satisfy \DE{} and \DW. By \cref{prop:dictionary}(b),(i),(f), $\bmu$ is finite, continuous from above, fulfils condition (E) and is weakly null-additive. By \cref{prop:A4}, \Lu$[\bmu]$ holds, which is \PL$[\mu]$ by \cref{prop:dictionary}(a).

%\section*{CRediT authorship contribution statement}
%[To be completed by the authors.]

%\section*{Declaration of competing interest}
%The authors declare that they have no known competing financial interests or personal relationships that could have appeared to influence the work reported in this paper.

%\section*{Data availability}
%No data was used for the research described in the article.

%\section*{Acknowledgements}

% AUTHORS: keep this statement only once you have actually done the verification it describes.
\section*{AI-use disclosure}
The authors used Claude (Opus 5.5) for language editing, organization, bibliographic checks, and checking calculations in arguments they had already developed. The authors reviewed all AI-assisted revisions, independently verified the mathematical arguments and proofs, and take full responsibility for the content of the paper.


\begin{thebibliography}{99}
\small

\bibitem{Bil68}
P.~Billingsley, Convergence of Probability Measures, Wiley, New York, 1968.

\bibitem{Cao15}
S.~Cao, Regularity of fuzzy measures on complete and separable metric spaces, in: 2015 12th International Conference on Fuzzy Systems and Knowledge Discovery (FSKD), IEEE, 2015, 168--172.

\bibitem{CM02}
A.~Castaldo, M.~Marinacci, A Lusin theorem for a class of Choquet capacities, Stat. Pap. 43(2002), 137--142.

\bibitem{Den94}
D.~Denneberg, Non-Additive Measure and Integral, Kluwer Academic Publishers, Dordrecht, 1994.

\bibitem{DF80}
I.~Dobrakov, J.~Farkov\'a, On submeasures II, Math. Slovaca 30(1980), 65--81.

\bibitem{HS73}
P.J.~Huber, V.~Strassen, Minimax tests and the Neyman--Pearson lemma for capacities, Ann. Stat. 1(1973), 251--263.

\bibitem{JS96}
Q.~Jiang, H.~Suzuki, Fuzzy measures on metric spaces, Fuzzy Sets Syst. 83(1996), 99--106.

\bibitem{JWZ99}
Q.~Jiang, S.~Wang, D.~Ziou, A further investigation for fuzzy measures on metric spaces, Fuzzy Sets Syst. 105(1999), 293--297.

\bibitem{Kec95}
A.S.~Kechris, Classical Descriptive Set Theory, Springer, New York, 1995.

\bibitem{Li03}
J.~Li, A further investigation for Egoroff's theorem with respect to monotone set functions, Kybernetika 39(2003), 753--760.

\bibitem{LM11}
J.~Li, R.~Mesiar, Lusin's theorem on monotone measure spaces, Fuzzy Sets Syst. 175(2011), 75--86.

\bibitem{LMPK15}
J.~Li, R.~Mesiar, E.~Pap, E.P.~Klement, Convergence theorems for monotone measures, Fuzzy Sets Syst. 281(2015), 103--127.

\bibitem{LOM21}
J.~Li, Y.~Ouyang, R.~Mesiar, Generalized convergence theorems for monotone measures, Fuzzy Sets Syst. 412(2021), 53--64.

\bibitem{LY04}
J.~Li, M.~Yasuda, Lusin's theorem on fuzzy measure spaces, Fuzzy Sets Syst. 146(2004), 121--133.

\bibitem{MUA04}
T.~Murofushi, K.~Uchino, S.~Asahina, Conditions for Egoroff's theorem in non-additive measure theory, Fuzzy Sets Syst. 146(2004), 135--146.

\bibitem{Pap95}
E.~Pap, Null-Additive Set Functions, Kluwer Academic Publishers, Dordrecht, 1995.

\bibitem{Sri98}
S.M.~Srivastava, A Course on Borel Sets, Springer, New York, 1998.

\bibitem{Sun94}
Q.~Sun, Property (S) of fuzzy measure and Riesz's theorem, Fuzzy Sets Syst. 62(1994), 117--119.

\bibitem{Tal08}
M.~Talagrand, Maharam's problem, Ann. Math. 168(2008), 981--1009.

\bibitem{Wal91}
P.~Walley, Statistical Reasoning with Imprecise Probabilities, Chapman and Hall, London, 1991.

\bibitem{WK09}
Z.~Wang, G.J.~Klir, Generalized Measure Theory, Springer, New York, 2009.

\bibitem{WT14}
T.~Watanabe, T.~Tanaka, On Lusin's theorem for non-additive measures that take values in an ordered topological vector space, Fuzzy Sets Syst. 244(2014), 41--50.

\bibitem{WYT13}
T.~Watanabe, S.~Yamada, T.~Tanaka, Note on the regularity of nonadditive measures, J. Appl. Math. 2013(2013), 234681.

\bibitem{Wie22}
J.~Wiesel, On a Lusin theorem for capacities, Proc. Am. Math. Soc. 150(2022), 1351--1360; arXiv:2010.09491.

\bibitem{WH94}
C.~Wu, M.~Ha, On the regularity of the fuzzy measure on metric fuzzy measure spaces, Fuzzy Sets Syst. 66(1994), 373--379.

\bibitem{YZH24}
D.~Yang, Z.~Zong, F.~Hu, Further results on Lusin's theorem for uncertain variables, Symmetry 16(2024), 1034.

\end{thebibliography}
\end{document}